\documentclass[11pt,reqno]{amsart}
\usepackage[T1]{fontenc}
\usepackage[utf8]{inputenc}
\usepackage{lmodern}
\usepackage[english]{babel}
\usepackage{amsmath,amssymb,amsfonts}
\usepackage[colorlinks=true,linkcolor=blue,citecolor=blue,urlcolor=blue]{hyperref}
\newtheorem*{theorem}{Theorem}
\newtheorem{lemma}{Lemma}[section]
\newtheorem{proposition}[lemma]{Proposition}
\newtheorem{corollary}[lemma]{Corollary}
\newcommand{\om}{\omega}
\newcommand{\Syl}{\operatorname{Syl}}
\newcommand{\Aut}{\operatorname{Aut}}
\newcommand{\cl}{\operatorname{cl}}
\newcommand{\rk}{\operatorname{rk}}
\newcommand{\F}{\mathbb F}
\title[\texorpdfstring{No solvable group is isospectral to $S_{10}$}{No solvable group is isospectral to S10}]{No solvable group is isospectral to the symmetric group of degree 10}
\author[D.V. Lytkina, V.D. Mazurov]{D.V. Lytkina, V.D. Mazurov}
\address{Daria Viktorovna Lytkina\newline\hphantom{iii}Novosibirsk State University,\newline\hphantom{iii}1 Pirogova St., Novosibirsk, 630090, Russia}
\email{daria.lytkin@gmail.com}
\address{Victor Danilovich Mazurov\newline\hphantom{iii}Sobolev Institute of Mathematics,\newline\hphantom{iii}4 Koptyug Ave., Novosibirsk, 630090, Russia}
\email{vic.mazurov@gmail.com}
\thanks{The work of D.V. Lytkina was supported by the Russian Science Foundation under grant no. 25-41-10008, \url{https://rscf.ru/en/project/25-41-10008/}. The work of V.D. Mazurov was carried out within the framework of the state assignment of the Sobolev Institute of Mathematics, project FWNF-2026-0017.}
\dedicatory{Dedicated to Lev Sergeevich Kazarin on his 80th birthday}
\date{}
\begin{document}
\sloppy
\begin{abstract}
We prove that no finite solvable group has the same spectrum as the symmetric group of degree~10. The proof analyzes the two possible Frobenius structures of a Hall $\{5,7\}$-subgroup and the actions of the corresponding subgroups on normal $\{2,3\}$-sections.
\end{abstract}
\maketitle

\section{Introduction and outline of the proof}
All groups considered in this paper are finite. The \emph{spectrum} of a group $X$ is the set of its element orders:
\[
\omega(X)=\{\,|x|\mid x\in X\,\}.
\]
The spectrum is closed under taking positive divisors and is uniquely determined by the subset $\mu(X)$ of its elements maximal under divisibility. Groups $X$ and $Y$ are called \emph{isospectral} if $\omega(X)=\omega(Y)$.

A group $X$ is \emph{recognizable by spectrum} if every finite group isospectral to $X$ is isomorphic to $X$. It is \emph{almost recognizable by spectrum} if there are finitely many isomorphism classes of finite groups with the same spectrum, and \emph{unrecognizable by spectrum} if there are infinitely many such classes.

The recognition problem asks to what extent the set of element orders determines a finite group up to isomorphism. Much of the research concerns simple groups and groups closely related to them. Results, methods, and open questions are surveyed by M.A. Grechkoseeva, V.D. Mazurov, W.J. Shi, A.V. Vasil'ev, and N.Y. Yang in~\cite{Survey}.

For symmetric groups $S_n$, substantial results were obtained by I.B. Gorshkov~\cite{Gorshkov} and subsequently by I.B. Gorshkov and A.N. Grishkov~\cite{GorshkovGrishkov}. In particular, recognizability was established in~\cite{Gorshkov} for
\[
n\notin\{2,3,4,5,6,8,10,15,16,18,21,27,33,35,39,45\},
\]
and in~\cite{GorshkovGrishkov} for $n\in\{15,16,18,21,27\}$. Together with the other results described in~\cite{Survey}, this gave the following picture: the groups $S_n$ are unrecognizable for
\[
n\in\{2,3,4,5,6,8\},
\]
and recognizable for all other $n$, with the possible exception of $n=10$. Thus $S_{10}$ was the only remaining case. It was singled out in~\cite[Problem~4.2]{Survey}: is the symmetric group $S_{10}$ recognizable by spectrum?

In~\cite{LytkinaMazurov}, we answered this question in the negative. We considered
\[
H=(S_9\times S_3)\cap A_{12}\cong(A_9\times C_3)\rtimes C_2,
\]
where $S_9$ and $S_3$ act on disjoint sets of nine and three letters, respectively, and proved that
\[
\omega(H)=\omega(S_{10}),\qquad
\mu(S_{10})=\{8,9,12,14,20,21,30\}.
\]
The group $H$ contains a nontrivial normal cyclic subgroup of order~3. By the Mazurov--Shi criterion~\cite{MazurovShi}, a finite group is unrecognizable by spectrum if and only if it is isospectral to a finite group with a nontrivial normal abelian subgroup. The construction therefore proves that $S_{10}$ is unrecognizable and completes the recognition problem for symmetric groups.

However, $H$ is nonsolvable, since it contains $A_9$. The existence of a solvable group with a prescribed spectrum is a question distinct from recognition: unrecognizability alone gives no information about the solvability of isospectral groups. This question for $S_{10}$ was posed separately in~\cite[Question~3.2]{LytkinaMazurov}. The present paper gives a negative answer to this question. Thus the spectrum of $S_{10}$ is realized by infinitely many pairwise nonisomorphic finite groups, but by no solvable group.

Our main result is as follows.
\begin{theorem}
There is no finite solvable group $G$ such that
\begin{equation}\label{spectrum}
\om(G)=\Omega:=\{1,2,3,4,5,6,7,8,9,10,12,14,15,20,21,30\}
=\om(S_{10}).
\end{equation}
\end{theorem}

The proof proceeds as follows. The absence of order $35$ forces a Hall $\{5,7\}$-subgroup to have one of two Frobenius structures: a $5$-kernel with complement of order $7$, or a $7$-kernel with complement of order $5$. In the first case we study the quotient by $O_{5'}(G)$, and in the second the quotient by $O_{7'}(G)$. In both cases, large elementary abelian subgroups control the normal $\{2,3\}$-part.

In the branch with a $7$-kernel, we exclude cyclic and generalized quaternion Sylow $2$-subgroups. The centralizer of an involution contains elements of orders $14$ and $30$. In the cyclic case, its normal odd-order complement has spectrum $\{1,3,5,7,15\}$. Its Frobenius or $2$-Frobenius structure yields order $28$ or $56$. In the quaternion case, a Hall $3'$-subgroup of this centralizer has a normal odd-order complement by Frobenius's normal $p$-complement theorem; the action of $Q_8$ again yields order $28$. Consequently, $G$ contains $C_2^2$.

A key tool common to both branches is a lemma on the action of an elementary abelian group of rank $4$ on a $3$-group: if the centralizer of every nonidentity automorphism has exponent at most $3$, then the whole $3$-group has exponent at most $3$. This excludes order $9$ from the normal $\{2,3\}$-part. An element of order $9$ then either retains its order in the quotient or has image of order $3$. For a $7$-kernel, the first possibility is excluded using the action on chief factors, and the second leads to two cases according to the normal $2$- and $3$-subgroups of the $\{2,3\}$-part.

For a $5$-kernel, the faithful action on the Frattini quotient and the Fong--Swan lift to characteristic zero allow us to apply a Shult-type theorem at the prime $7$. The resulting normal Sylow $7$-subgroup in the complement is used to analyze its centralizer. If the centralizer has odd order, the $2$-part of the complement has order at most $2$. If it has even order, the absence of elements of order $6$ restricts its structure. Automorphisms of order $4$ and subgroups isomorphic to $C_2^2$ are then excluded separately.

The special external results used in these reductions are stated in Subsection~\ref{sec:external}.

\section{Notation and standard tools}
Assume that $G$ is solvable and satisfies~\eqref{spectrum}. Write $\pi(X)$ for the set of prime divisors of $|X|$ and $|X|_p$ for the largest power of $p$ dividing $|X|$. The subgroup $O_\pi(X)$ is the largest normal $\pi$-subgroup of $X$; $O_{p'}(X)$ corresponds to all primes other than $p$. A Hall $\pi$-subgroup has $\pi$-order and $\pi'$-index. The Fitting subgroup $F(X)$ is the largest normal nilpotent subgroup; for a finite group, $F(X)=\prod_pO_p(X)$.

We write $\Phi(P)$ for the intersection of the maximal subgroups of $P$, and $\gamma_1(P)=P$, $\gamma_{i+1}(P)=[\gamma_i(P),P]$ for its lower central series. The centralizer and normalizer of $Y$ in $X$ are denoted by $C_X(Y)$ and $N_X(Y)$. The notation $C_p^d$ denotes an elementary abelian group of order $p^d$ and rank $d$; $X^\#=X\setminus\{1\}$. The exponent $\exp X$ is the least common multiple of its element orders. An action is faithful if its kernel is trivial, and fixed-point-free if its centralizer in the group is trivial (or its fixed-point subspace in a module is zero).

Every element order in a subgroup or quotient of $G$ belongs to $\Omega$. For a quotient, this follows because the order of an image divides the order of a preimage and $\Omega$ is closed under divisibility. In particular, the following orders are forbidden:
\[
18,24,28,35,40,42,45,56,60,63.
\]
A series $1=G_0\trianglelefteq G_1\trianglelefteq\cdots\trianglelefteq G_t=G$ is a $G$-chief series if every $G_i$ is normal in $G$ and there is no intermediate normal subgroup of $G$ between consecutive terms. For solvable $G$, its factors $G_i/G_{i-1}$ are elementary abelian. A chief factor below $R$ means a factor of such a series lying between $1$ and the normal subgroup $R$.

We use $x^g=g^{-1}xg$ and $[x,y]=x^{-1}y^{-1}xy$, with longer commutators left-normed. The subgroup $[X,A]$ is generated by $x^{-1}x^a$ for $x\in X$, $a\in A$. An action of $A$ on $X$ is \emph{coprime} if it is by automorphisms and $(|A|,|X|)=1$. Unless otherwise stated, all groups considered below are solvable.

\subsection{External results}\label{sec:external}
We use standard results on finite groups and their representations in their usual forms, giving references where they are applied. Here we collect the special results needed for the structural reductions.

\begin{proposition}[Gruenberg--Kegel, {\cite{Williams}}]\label{extGK}
Let $X$ be a finite solvable group with disconnected prime graph: its vertices are the prime divisors of $|X|$, and $r,s$ are adjacent if and only if $rs\in\om(X)$. Then $X$ is Frobenius or $2$-Frobenius. In the latter case there exist $1<K<L<X$, with $K,L\trianglelefteq X$, such that $L$ is Frobenius with kernel $K$ and $X/K$ is Frobenius with kernel $L/K$. In the Frobenius case, the components are $\pi(K)$ and $\pi(H)$, where $K$ is the kernel and $H$ a complement.
\end{proposition}

\begin{proposition}[Flavell, {\cite[Theorem A]{Flavell}}]\label{extFlavell}
Let $r$ be a prime, let a group $A$ of order $r$ act on a finite $r'$-group $Y$, and let $M$ be a faithful irreducible module for $Y\rtimes A$ over a field of arbitrary characteristic. If $C_M(A)=0$, then either $[Y,A]=1$, or $r$ is a Fermat prime and $[Y,A]$ is a nonabelian special $2$-group. Here a $2$-group $S$ is special if $S'=\Phi(S)=Z(S)$ and this common subgroup is elementary abelian.
\end{proposition}

\begin{proposition}[Espuelas, {\cite[p.~207]{Espuelas}}]\label{extEspuelas}
Let $X$ be a finite $p$-solvable group and let $M$ be a faithful $KX$-module with $\operatorname{char}K\nmid|O_{p'}(X)|$. Suppose that $x\in X$ has order $p^n$, $n\ge1$, and $C_M(x)=0$. If $p$ is a Fermat prime, assume also that the Sylow $2$-subgroups of $X$ are abelian. If $p=2$, assume that the Sylow $q$-subgroups of $X$ are abelian for every Mersenne prime $q<2^n$. Then
\[
x^{p^{n-1}}\in O_p(X).
\]
A Fermat prime has the form $2^{2^a}+1$, and a Mersenne prime has the form $2^a-1$.
\end{proposition}

\subsection{Elementary consequences and technical tools}
Over a splitting field $k$ with $\operatorname{char}k\nmid |E|$, a module for an elementary abelian group $E$ is a direct sum of weight spaces for characters $\lambda:E\to k^\times$. Such a one-dimensional character is called a \emph{weight}. For an arbitrary group, an isotypic component is the sum of all simple submodules of a fixed isomorphism type. If $E\cong C_\ell^d$, the kernel of a weight has codimension at most $1$, since a finite subgroup of $k^\times$ is cyclic. Extension of scalars does not change the dimension of the kernel of a linear operator: it equals the dimension of the space minus the rank of its matrix.

\begin{lemma}
\label{fielddegree}
Let $p\ne q$ be primes and let $C_q$ act fixed-point-freely on a nonzero elementary abelian $p$-module $M$. Then $\operatorname{ord}_q(p)$ divides $\dim_{\F_p}M$. If an elementary abelian $q$-group $E$ acts on a nonzero module of characteristic different from $q$, and every $e\in E^\#$ acts fixed-point-freely, then $|E|\le q$.
\end{lemma}
\begin{proof}
For the first assertion, extend scalars to an algebraic closure of $\F_p$. The eigenvalues of a generator of $C_q$ are $q$th roots of unity, and $1$ does not occur. All the others have order $q$. The Frobenius automorphism $\lambda\mapsto\lambda^p$ partitions them into orbits of length $d=\operatorname{ord}_q(p)$, since $\lambda^{p^i}=\lambda$ is equivalent to $p^i\equiv1\pmod q$. The characteristic polynomial is defined over $\F_p$, so roots in the same orbit have equal multiplicities. The dimension is therefore a sum of multiples of $d$.

For the second assertion, choose a weight $\lambda:E\to k^\times$ occurring over a splitting field. Its image has order $1$ or $q$, so $|E:\ker\lambda|\le q$. If $|E|>q$, its kernel contains $e\ne1$, which fixes a nonzero vector in the weight space, a contradiction.
\end{proof}

\begin{lemma}
\label{inversionlemma}
If an automorphism $a$ of a finite group $X$ has order $2$ and $C_X(a)=1$, then $X$ is abelian of odd order and $x^a=x^{-1}$ for every $x\in X$.
\end{lemma}
\begin{proof}
The map $f:X\to X$, $f(u)=u^{-1}u^a$, is injective: $f(u)=f(v)$ implies $vu^{-1}=(vu^{-1})^a$, whence $u=v$. Since $X$ is finite, $f$ is surjective. If $x=f(u)$, then $x^a=(u^a)^{-1}u=x^{-1}$. Thus inversion is an automorphism, and for $x,y\in X$ we obtain $x^{-1}y^{-1}=(xy)^{-1}=y^{-1}x^{-1}$. Hence $X$ is abelian. It has no involutions, since an involution would be fixed by inversion; Cauchy's theorem shows that $|X|$ is odd.
\end{proof}

\begin{lemma}
\label{Hallsections}
Let $N\trianglelefteq X$, let $X$ be solvable, and let $H$ be a Hall $\pi$-subgroup of $X$. Then $H\cap N$ is a Hall $\pi$-subgroup of $N$, and $HN/N$ is a Hall $\pi$-subgroup of $X/N$. If $N$ is a $p'$-group, the map from any Sylow $p$-subgroup of $X$ to its image is an isomorphism. If $X/N$ is a $p$-group, every Sylow $p$-subgroup of $X$ maps onto $X/N$.
\end{lemma}
\begin{proof}
The equality $|N:H\cap N|=|HN:H|$ shows that this index divides $|X:H|$ and is a $\pi'$-number. Hence $H\cap N$ is a Hall subgroup of $N$. Similarly, $|X/N:HN/N|=|X:HN|$ is a $\pi'$-number, whereas $HN/N$ is a $\pi$-group. The remaining assertions follow from $|S N/N|=|S|/|S\cap N|$ by comparing $p$-parts.
\end{proof}

\begin{lemma}
\label{normalizerlift}
Let $R\trianglelefteq X$ be a $p'$-subgroup of a solvable group, and let $\bar P\trianglelefteq X/R$ be a $p$-subgroup. Its full preimage has a complement $P$ to $R$, and
\[
X=R N_X(P),\quad P\cong\bar P,\quad
\ker(N_X(P)\to X/R)=C_R(P).
\]
\end{lemma}
\begin{proof}
In the preimage of $\bar P$, the Schur--Zassenhaus theorem~\cite[9.1.2]{Robinson} gives $P$ and conjugacy of all complements by elements of $R$. For $x\in X$, the subgroup $P^x$ is another such complement. Hence $P^x=P^r$ for some $r\in R$, so $xr^{-1}\in N_X(P)$ and $X=R N_X(P)$. Finally, if $r\in R\cap N_X(P)$, then $[P,r]\le P\cap R=1$, and $r\in C_R(P)$. The reverse inclusion is immediate.
\end{proof}

\begin{lemma}
\label{localcentralizer}
Let $A$ be a group of automorphisms of a finite group $X$, and let $N\trianglelefteq X$ be $A$-invariant with $(|A|,|N|)=1$. Then
\[
C_{X/N}(A)=C_X(A)N/N.
\]
\end{lemma}
\begin{proof}
The right-to-left inclusion follows directly from the definition. Let $xN$ be fixed by $A$. Consider $X\rtimes A$ with an external copy of $A$. The condition $[x,A]\le N$ means that $x$ normalizes $NA$. Both $A$ and $A^x$ complement $N$ in $NA$, since both have order $|A|$ and intersect $N$ trivially. By Schur--Zassenhaus~\cite[9.1.2]{Robinson}, there is $n\in N$ with $A^x=A^n$. Then $y=xn^{-1}$ normalizes $A$. For $a\in A$, the commutator $[y,a]$ lies in $A$ by normalization and in $N$ because $yN=xN\in C_{X/N}(A)$. Thus $[y,A]=1$, since $A\cap N=1$. We have $y\in C_X(A)$ and $yN=xN$.
\end{proof}

\begin{lemma}
\label{norm}
Let $P$ be a normal $p$-subgroup of $P\langle c\rangle$, and let $M$ be a module of characteristic $r\ne p$ over a splitting field. If $c$ cyclically permutes $m$ distinct $P$-isotypic components, then the operator
\[
1+c+\cdots+c^{m-1}
\]
is nonzero on their sum. If $c^m$ acts trivially on this sum, then $c$ has a nonzero fixed point. In the group-theoretic application, a nonzero norm operator on an elementary abelian $r$-section $L/K$ normalized by $c$ yields an element whose order is divisible by $rm$, provided that $|c|=m$ and the image of $c$ modulo $L$ still has order $m$.
\end{lemma}
\begin{proof}
Denote the components by $M_0,\ldots,M_{m-1}$, where $cM_i=M_{i+1}$ with indices modulo $m$. For $0\ne v\in M_0$, the vectors $v,cv,\ldots,c^{m-1}v$ lie in distinct direct summands. Thus their sum $w$ is nonzero: its projection onto $M_0$ is $v$. If $c^m$ acts trivially on this sum, then
\[
cw=cv+c^2v+\cdots+c^mv=w.
\]

For the group-theoretic application, let $M=L/K$ be an elementary abelian $r$-section normalized by $c$ of order $m$, and suppose that $c$ retains order $m$ modulo $L$. Choose $v\in L$ whose image $\bar v$ has nonzero norm. Multiplication in the extension gives
\[
(vc)^mK=\bar v+c\bar v+\cdots+c^{m-1}\bar v\in L/K.
\]
The right-hand side has order $r$. If $n=|vc|$, then $m\mid n$ because the image of $vc$ modulo $L$ has order $m$. Since $(vc)^m$ has image of order $r$, we have $r\mid n/m$. Hence $rm\mid n$, and $\langle vc\rangle$ contains an element of order $rm$.
\end{proof}

\begin{lemma}
\label{characters}
Let an automorphism $a$ of prime order $q\ne p$ act fixed-point-freely on a $p$-group $P$. Then the principal character is the only irreducible character of $P$ fixed by $a$. The same holds for irreducible types over a splitting field of characteristic different from $p$.
If two commuting automorphisms of distinct prime orders $q,s$ act fixed-point-freely on $P$, all nonprincipal types have orbits of length $qs$ under their product.
\end{lemma}
\begin{proof}
Let $\mathcal K$ be a conjugacy class of $P$ invariant under $a$. Its size is a power of $p$ and is not divisible by $q$. The orbits of $\langle a\rangle$ on $\mathcal K$ have length $1$ or $q$. If there were no fixed elements, then $q$ would divide $|\mathcal K|$. Hence $\mathcal K$ contains $u\in C_P(a)$. By hypothesis, $u=1$, and $\mathcal K=\{1\}$.

Thus $a$ fixes exactly one conjugacy class of $P$. By Brauer's permutation lemma~\cite[Theorem 6.32]{Isaacs}, it fixes exactly one irreducible character. The principal character is fixed, so there are no others. Passage to a splitting field of characteristic $r\ne p$ is justified by the coprime correspondence between ordinary and Brauer characters~\cite[Theorem 15.13]{Isaacs}, applied to the $p$-group $P$: $r$ does not divide $|P|$, and the correspondence of types is compatible with automorphisms.

For the last assertion, let $a,b$ commute and have distinct prime orders $q,s$. Then $\langle a,b\rangle$ is cyclic of order $qs$. By the first part, the stabilizer of a nonprincipal type cannot contain a subgroup of order $q$ or $s$. Every nontrivial subgroup of a cyclic group of order $qs$ contains one of these subgroups. Thus the stabilizer is trivial and the orbit has length $qs$.
\end{proof}

\begin{lemma}
\label{commutator}
Let a group $A$ act coprimely on a finite $\{2,3\}$-group $R$, and let $1=R_0<\cdots<R_t=R$ be a series whose terms are normal in $R\rtimes A$ and whose factors are elementary abelian. If $A$ centralizes all $3$-factors of this series, then
\[
[R,A]\le O_2(R).
\]
Interchanging $2$ and $3$ gives the dual assertion. If $A$ centralizes all factors of the series, then $[R,A]=1$.
\end{lemma}
\begin{proof}
Refine the given series so that each factor is a minimal normal subgroup for the action of $R\rtimes A$; centralization of its $3$-factors is preserved. Argue by induction on $|R|$. The assertion is clear for $R=1$. Let $L$ be the first nontrivial term of the refined series. Solvability and minimality imply that $L$ is elementary abelian, and is a $2$-group or a $3$-group.

Set $D=[R,A]$. This subgroup is normal in $R$ and $A$-invariant, and $[D,A]=D$ by coprime action~\cite[8.2.7(b)]{KurzweilStellmacher}. Induction shows that $DL/L$ is a $2$-group. If $L$ is a $2$-group, both $D\cap L$ and $D/(D\cap L)$ are $2$-groups, so $D$ is a normal $2$-subgroup of $R$.

Suppose that $L$ is a $3$-group. The hypothesis gives $[L,A]=1$. For $r\in R$ and $a\in A$, conjugation by $r$ and by $r^a$ induces the same action on $L$, since $L\trianglelefteq R$ and $a$ fixes every element of $L$. Thus $[R,A]$ centralizes $L$. Consequently, $D\cap L$ is a central normal Sylow $3$-subgroup of $D$, and $D/(D\cap L)$ is a $2$-group.

By Schur--Zassenhaus~\cite[9.1.2]{Robinson}, $D\cap L$ has a complement $D_2$. Since $D\cap L$ is central, $D=(D\cap L)\times D_2$, and $D_2$ is the unique Sylow $2$-subgroup of $D$, hence is $A$-invariant. The group $A$ centralizes $D\cap L$ and preserves $D_2$, so $[D,A]\le D_2$. But $[D,A]=D$, whence $D\cap L=1$, and $D$ is again a $2$-group.

In both cases $D\le O_2(R)$. Interchanging the primes gives the dual assertion. If all factors are centralized, $D$ lies in both $O_2(R)$ and $O_3(R)$, whose intersection is trivial.
\end{proof}

\section{The rank-four lemma}
We first bound the nilpotency class and then the exponent of a $3$-group under an elementary abelian group of automorphisms.

\begin{lemma}
\label{rankfour}
Let $E\cong C_\ell^4$, where $\ell\ne3$ is prime, act on a finite $3$-group $T$. If
\[
\exp C_T(e)\le3\qquad(e\in E^\#),
\]
then $\exp T\le3$.
\end{lemma}
\begin{proof}
By the classical result of Levi~\cite[12.3.5--12.3.6]{Robinson}, every group of exponent $3$ has class at most $3$ and satisfies the identities $[u,v,v]=1$ and $[u,v,w][u,w,v]=1$. We apply these identities to $C_T(e)$ for $e\in E^\#$.

Consider the graded Lie algebra over $\F_3$
\[
L=\bigoplus_{i\ge1}\gamma_i(T)/(\gamma_{i+1}(T)\gamma_i(T)^3).
\]
Its bracket is induced by the group commutator: $[\overline u,\overline v]=\overline{[u,v]}$. The inclusions $[\gamma_i(T),\gamma_j(T)]\le\gamma_{i+j}(T)$ ensure compatibility of degrees, and factoring out cubes makes the components vector spaces over $\F_3$. The usual commutator identities modulo the next lower central term give bilinearity and the Jacobi identity. Since $\gamma_i(T)/\gamma_{i+1}(T)$ is generated by weight-$i$ commutators in generators of $T$, the algebra $L$ is generated by its degree-one component. Extend scalars to a splitting field for $E$. Any three weights of $E$ have a nontrivial common kernel: each kernel has codimension at most $1$, whereas $\rk E=4$.

For $e\in E^\#$, coprime lifting of fixed points gives
\[
C_{U/W}(e)=C_U(e)W/W
\]
for every $e$-invariant section $U/W$ of $T$. Thus the fixed subspace of each homogeneous component of $L$ is the image of $C_{\gamma_i(T)}(e)$. Levi's two group identities induce
\[
[x,y,y]=0,\qquad [x,y,z]+[x,z,y]=0
\]
on these homogeneous fixed subspaces. The second identity is trilinear and therefore remains valid after extending scalars. On a fixed homogeneous subspace the first identity extends by polarization, since $2\ne0$ in characteristic $3$.

Now let $x,y,z$ be homogeneous in degree and weight. The three weights have a common nonidentity kernel element $e\in E$, so the preceding identities apply to $x,y,z$. Expanding arbitrary elements of $L$ into homogeneous weight vectors gives $[x,y,y]=0$ throughout $L$: diagonal terms vanish by the first identity, and the remaining terms cancel in pairs by the second.

Set $X=\operatorname{ad}x$, $Y=\operatorname{ad}y$, $Z=\operatorname{ad}z$. The polarized identity gives $XY=-YX$, so $\operatorname{ad}[x,y]=[X,Y]=2XY$. Applying the same anticommutation to $\operatorname{ad}[x,y]$ and $Z$, we obtain
\[
0=2XYZ+2ZXY=4XYZ.
\]
In characteristic $3$, this implies $XYZ=0$. Thus $\cl L\le3$ and $L_4=0$. For the finite abelian $3$-group $W_4=\gamma_4(T)/\gamma_5(T)$, it follows that $W_4=3W_4$, hence $W_4=0$. Indeed, if $|W_4|=3^b$, iteration yields $W_4=3^bW_4=0$. The equality $\gamma_4(T)=\gamma_5(T)$ implies $\gamma_i(T)=\gamma_4(T)$ for every $i\ge4$. Since a finite $3$-group is nilpotent, its lower central series terminates at the identity, and $\gamma_4(T)=1$.

We now establish the exponent. Put
\[
\mathcal X=\bigcup_{B<E,\ |E:B|=\ell}C_T(B).
\]
By~\cite[8.3.4(a)]{KurzweilStellmacher}, $T$ is generated by the centralizers $C_T(B)$ for $B\le E$ with $E/B$ cyclic. Since $E$ is elementary abelian, such a $B$ is either a hyperplane or $E$ itself; moreover, $C_T(E)$ is contained in every hyperplane centralizer. Hence $T=\langle\mathcal X\rangle$. If $x_i\in C_T(B_i)$ for $i=1,2,3$, then $\dim(B_1\cap B_2\cap B_3)\ge4-3=1$. Thus any three elements of $\mathcal X$ are simultaneously centralized by some $e\in E^\#$ and generate a group of exponent at most $3$.

Since $\cl T\le3$, the group $T'$ is abelian and generated by commutators of weights $2$ and $3$ in elements of $\mathcal X$. Each has order dividing $3$, so $\exp T'\le3$. Commutators of weight $3$ are central. The identities $[uv,w]=[u,w]^v[v,w]$ and $[u,vw]=[u,w][u,v]^w$ therefore show that a triple commutator is trilinear in the images of its arguments in $T/T'$. The identities $[x,y,y]=1$ and $[x,y,z][x,z,y]=1$, established on $\mathcal X$, consequently extend to all of $T$.

Thus $T$ is a $2$-Engel group of class at most $3$. For any $u,v$, the equalities $[u,v,v]=1$ and $[v,u,u]=1$ mean that $[u,v]$ commutes with both $u$ and $v$, so $\langle u,v\rangle$ has class at most $2$. The cube formula and $\exp T'\le3$ give
\[
(uv)^3=u^3v^3[v,u]^3=u^3v^3.
\]
The cube map is a homomorphism and is trivial on $\mathcal X$, hence on $T$.
\end{proof}

\begin{lemma}
\label{twoelementary}
Let $R\trianglelefteq G$ be a $\{2,3\}$-group, and let $E\cong C_7^3$ be a subgroup of $G$. Then the Sylow $2$-subgroups of $R$ are elementary abelian.
\end{lemma}
\begin{proof}
Since $R\trianglelefteq G$, the group $E$ acts coprimely on $R$. By coprime action~\cite[8.2.3(a)]{KurzweilStellmacher}, choose an $E$-invariant $S\in\Syl_2(R)$. If $e\in E^\#$ and $u\in C_S(e)$ had order $4$, then $ue$ would have order $28$. Thus every element of $C_S(e)$ has square $1$. A group of exponent $2$ is abelian, since $(uv)^2=1$ implies $uv=(uv)^{-1}=v^{-1}u^{-1}=vu$.

By coprime action~\cite[8.3.4(a)]{KurzweilStellmacher}, as in the proof of Lemma~\ref{rankfour}, $S$ is generated by the union $\mathcal X$ of the centralizers $C_S(B)$ over hyperplanes $B<E$. Each $u\in\mathcal X$ is an involution or the identity: $B$ has rank $2$ and contains $e\ne1$, so $u\in C_S(e)$.

For $u\in C_S(B_1)$ and $v\in C_S(B_2)$, the intersection $B_1\cap B_2$ has rank at least $1$, since $\rk E=3$. Choose $e\ne1$ in it. Then $u,v\in C_S(e)$ and $[u,v]=1$. Hence $S$ is generated by pairwise commuting elements of order at most $2$ and is elementary abelian. All Sylow $2$-subgroups of $R$ are conjugate to $S$ and have the same property.
\end{proof}

\section{\texorpdfstring{A Hall $\{5,7\}$-subgroup and preliminary reductions}{A Hall {5,7}-subgroup and preliminary reductions}}
\begin{lemma}
\label{orientations}
A Hall $\{5,7\}$-subgroup $H$ has one of the following forms:
\begin{enumerate}
\item[(A)] $H=P\rtimes C_7$, where $P\in\Syl_5(G)$, $\exp P=5$, and $\rk Z(P)$ is a positive multiple of $6$;
\item[(B)] $H=Q\rtimes C_5$, where $Q\in\Syl_7(G)$, $\exp Q=7$, and $\rk Z(Q)$ is a positive multiple of $4$.
\end{enumerate}
In case \textup{(A)}, a Sylow $7$-subgroup of $G$ is cyclic of order $7$; in case \textup{(B)}, a Sylow $5$-subgroup is cyclic of order $5$.
\end{lemma}
\begin{proof}
Hall's theorem~\cite[9.1.7]{Robinson} gives a Hall $\{5,7\}$-subgroup $H$ of the solvable group $G$. Both primes divide $|H|$, and $35\notin\om(H)$. Its prime graph is therefore disconnected, so Proposition~\ref{extGK} gives the Frobenius or $2$-Frobenius case.

Suppose that $H$ is $2$-Frobenius. By the standard properties of such groups~\cite[Lemmas 2.1--2.5]{CostanzoLewis}, its middle complement $A$ is cyclic, of odd order, and a Hall subgroup, while the nontrivial upper complement $C$ acts faithfully on $A$. Since orders $25,49,35$ are absent, $|A|=5$ or $7$. As $A$ is Hall, $C$ is a $7$-group in the first case and a $5$-group in the second. However, $|C|$ must divide $|\Aut(C_5)|=4$ or $|\Aut(C_7)|=6$, respectively. Neither allows $C\ne1$, a contradiction.

Thus $H=K\rtimes A$ is Frobenius. Thompson's theorem on Frobenius kernels~\cite[10.5.6]{Robinson} shows that $K$ is nilpotent. It cannot contain both primes $5,7$, since commuting Sylow subgroups would yield order $35$. The kernel and complement have coprime orders: if $r\mid|A|$ is prime, an element of order $r$ acts fixed-point-freely on $K$, and counting its orbits gives $|K|\equiv1\pmod r$. Hence $K$ is a full Sylow subgroup for one prime and $A$ for the other.

Odd-order Sylow subgroups of a Frobenius complement are cyclic~\cite[8.5.4]{Robinson}. The absence of orders $25,49$ therefore gives $|A|=7$ for a $5$-kernel and $|A|=5$ for a $7$-kernel. The kernel has exponent $5$ or $7$, since larger powers of that prime are forbidden.

The center of the nontrivial kernel is a nonzero elementary abelian group invariant under the complement, which acts fixed-point-freely on it. Lemma~\ref{fielddegree} gives a dimension divisible by $\operatorname{ord}_7(5)=6$ in the first case and by $\operatorname{ord}_5(7)=4$ in the second. Finally, $H$ is Hall in $G$, so its Sylow subgroups are also Sylow subgroups of $G$. This proves all assertions in (A) and (B).
\end{proof}

\begin{lemma}
\label{radical}
In orientation \textup{(A)}, $R=O_{5'}(G)$ is a $\{2,3\}$-group; in orientation \textup{(B)}, the same holds for $R=O_{7'}(G)$. For $B=G/R$, respectively,
\[
F(B)=O_5(B)\ne1\quad\text{or}\quad F(B)=O_7(B)\ne1,
\qquad C_B(F(B))\le F(B).
\]
\end{lemma}
\begin{proof}
First consider (A). Let $R=O_{5'}(G)$, and let $H=P_5\rtimes Q_7$ be a Hall $\{5,7\}$-subgroup. If $7\mid|R|$, Lemma~\ref{Hallsections} shows that the normal intersection $R\cap H$ contains a full Sylow $7$-subgroup of $H$, of order $7$. All complements of order $7$ in $H$ are conjugate, and $R\cap H$ is normal in $H$. Thus we may assume $Q_7\le R\cap H$. Then $[P_5,Q_7]\le R\cap H$.

The coprime action of $Q_7$ on $P_5$ is fixed-point-free. Coprime action~\cite[8.2.7(a)]{KurzweilStellmacher} gives $P_5=C_{P_5}(Q_7)[P_5,Q_7]=[P_5,Q_7]$. Hence $P_5\le R$, contrary to $R$ being a $5'$-group. Thus $7\nmid|R|$ and $\pi(R)\subseteq\{2,3\}$.

In (B), interchange $5$ and $7$. If $5\mid|O_{7'}(G)|$, its normal intersection with $Q_7\rtimes C_5$ contains a complement $C_5$ and hence $[Q_7,C_5]=Q_7$, again a contradiction.

Now let $p$ be the kernel prime, $p=5$ or $7$, and let $B=G/O_{p'}(G)$. If $O_{p'}(B)\ne1$, its full preimage would be a normal $p'$-subgroup of $G$ properly containing $O_{p'}(G)$. Thus $O_{p'}(B)=1$. All normal $q$-subgroups of $B$ for $q\ne p$ are trivial, and $F(B)=\prod_qO_q(B)=O_p(B)$. Standard properties of the Fitting subgroup~\cite[5.4.4]{Robinson} show that this subgroup is nontrivial, since $B$ is nontrivial and solvable, and give $C_B(F(B))\le F(B)$.
\end{proof}

\begin{lemma}
\label{noradicalnine}
For $R$ in Lemma~\ref{radical}, we have $9\notin\om(R)$.
\end{lemma}
\begin{proof}
In (A), the center of a Sylow $5$-subgroup contains $E\cong C_5^4$, since its rank is at least $6$. In (B), the center of a Sylow $7$-subgroup contains $E\cong C_7^4$. Normality of $R$ gives a conjugation action of $E$ on $R$ in both cases. By Lemma~\ref{radical}, $R$ is a $\{2,3\}$-group, so the action is coprime.

Choose an $E$-invariant $T\in\Syl_3(R)$ by coprime action~\cite[8.2.3(a)]{KurzweilStellmacher}. For $e\in E^\#$, the centralizer $C_T(e)$ has no element of order $9$: commuting elements of orders $9$ and $5$ in (A) would yield $45$, and orders $9$ and $7$ in (B) would yield $63$. Since all element orders in $T$ are powers of $3$ and $27\notin\Omega$, we have $\exp C_T(e)\le3$.

The rank of $E$ is $4$, so Lemma~\ref{rankfour} gives $\exp T\le3$. Every $3$-element of $R$ is conjugate in $R$ to an element of $T$. Thus $R$ has no element of order $9$.
\end{proof}

\section{\texorpdfstring{The branch with a $7$-kernel}{The branch with a 7-kernel}}
\subsection{\texorpdfstring{Excluding two types of Sylow $2$-subgroups}{Excluding two types of Sylow 2-subgroups}}
The final part of branch (B) requires a subgroup $C_2^2$ in $G$. The next lemma isolates the required information about the odd-order complement.

\begin{lemma}
\label{oddclassification}
Let $N$ be a solvable group of odd order with $\om(N)=\{1,3,5,7,15\}$. Then either a Sylow $7$-subgroup of $N$ has order $7$, or
\[
N=Q_0\rtimes C_{15},
\]
where $Q_0$ is a normal Sylow $7$-subgroup and $C_{15}$ acts fixed-point-freely on it.
\end{lemma}
\begin{proof}
The prime graph of $N$ has components $\{3,5\}$ and $\{7\}$. By Proposition~\ref{extGK}, $N$ is Frobenius or $2$-Frobenius.

In the Frobenius case, the prime divisors of the kernel and complement form the two components. If the complement is a $7$-group, it is cyclic, since odd-order Sylow subgroups of a Frobenius complement are cyclic. The absence of $49$ then gives order $7$, and it is a full Sylow $7$-subgroup of $N$. If the kernel is a $7$-group, the complement has prime divisors $3,5$. Its Sylow subgroups are cyclic, and the absence of $9,25$ gives order $15$. A group of order $15$ is cyclic: by Sylow's theorems, the number of Sylow $5$-subgroups divides $3$ and is congruent to $1$ modulo $5$, and the number of Sylow $3$-subgroups divides $5$ and is congruent to $1$ modulo $3$. Both numbers are $1$, so the group is $C_5\times C_3$.

Suppose that $N$ is $2$-Frobenius, with $1<K<L<N$, where $L$ is Frobenius with kernel $K$ and $N/K$ is Frobenius with kernel $L/K$. The middle factor $L/K$ is isomorphic to a complement $H$ in $L$, which is also the Frobenius kernel of its normalizer. Hence $H$ is nilpotent and its odd-order Sylow subgroups are cyclic. Thus $H$ is cyclic and is a Hall subgroup of $N$. These standard properties are collected in~\cite[Lemmas 2.1--2.5]{CostanzoLewis}.

The possible orders of $H$ are $3$, $5$, $7$, and $15$. The nontrivial upper complement $N/L$ has odd order and embeds in $\Aut(H)$. For $|H|=3,5,15$, the group $\Aut(H)$ is a $2$-group. Therefore $|H|=7$, and since $H$ is Hall, it is a full Sylow $7$-subgroup of $N$.
\end{proof}

\begin{lemma}
\label{nocyclictwo}
A Sylow $2$-subgroup of $G$ cannot be cyclic.
\end{lemma}
\begin{proof}
Suppose that $S=\langle g\rangle\in\Syl_2(G)$ is cyclic. Since $8\in\Omega$ and $16\notin\Omega$, we have $|S|=8$. Put $z=g^4$ and $C=C_G(z)$. All involutions of $G$ are conjugate to $z$: each lies in a conjugate of $S$, which has a unique involution. Conjugating elements of orders $14$ and $30$ so that their involutions become $z$, we obtain $14,30\in\om(C)$.

The group $S$ is also a Sylow $2$-subgroup of $C$. Burnside's normal $p$-complement criterion~\cite[10.1.8]{Robinson} gives a normal odd-order complement $N$ in $C$: the odd-order group $N_C(S)/C_C(S)$ embeds in the $2$-group $\Aut(S)$ and is therefore trivial. Thus $C=N\rtimes S$.

Since $z\in Z(C)$, orders $9,21$ do not occur in $N$, for otherwise $18,42$ would occur together with $z$. The presence of $14,30$ in $C$ gives $7,15\in\om(N)$. Hence
\[
\om(N)=\{1,3,5,7,15\}.
\]
Apply Lemma~\ref{oddclassification}.

If a Sylow $7$-subgroup of $N$ has order $7$, coprime action gives an $S$-invariant such subgroup $Q_0$. The action of $S$ on $Q_0$ has order at most $2$, since $\Aut(C_7)\cong C_6$. Thus $g^2$, of order $4$, centralizes $Q_0$, yielding order $28$.

Let $N=Q_0\rtimes A$, where $A\cong C_{15}$. Choose a Hall $\{2,3,5\}$-subgroup $J$ of $C$ containing $S$. Then $A_1=J\cap N$ has order $15$, is normal in $J$, and is conjugate in $N$ to $A$. Hence $A_1\cong C_{15}$ acts fixed-point-freely on $Q_0$ and is normalized by $g$. Denote its subgroups of orders $3,5$ by $A_3,A_5$.

The element $g$ inverts $A_3$; otherwise it centralizes $A_3$ and yields order $24$. Thus $g^2$ centralizes $A_3$. The action of $g$ on $A_5$ must have order $4$: if it had order at most $2$, then $g^2$ would centralize all of $A_1$, and an element of order $4$ together with $A_1$ would yield order $60$.

The subgroup $Q_0$ is normal in $C$, being the unique Sylow $7$-subgroup of the normal subgroup $N$. The group $A_5$ acts fixed-point-freely on the nonzero module $M=Q_0/\Phi(Q_0)$. Over a splitting field, $g$ permutes its nontrivial weights in orbits of length $4$. Moreover, $g^4=z$ acts trivially on $Q_0$, since $Q_0\le C$. Applying Lemma~\ref{norm} to $g$ on $M$, with $g^4=1$, gives $C_M(g)\ne0$. Coprime lifting yields $C_{Q_0}(g)\ne1$. An element of order $7$ in this centralizer, together with $g$, yields order $56$.
\end{proof}

\begin{lemma}
\label{noquaterniontwo}
A Sylow $2$-subgroup of $G$ cannot be generalized quaternion.
\end{lemma}
\begin{proof}
Suppose that $S\in\Syl_2(G)$ is generalized quaternion. The presence of order $8$ and absence of $16$ give $S\cong Q_{16}$. Let $z$ be its unique involution. As in Lemma~\ref{nocyclictwo}, all involutions of $G$ are conjugate to $z$, and $C=C_G(z)$ contains elements of orders $14,30$.

Take a Hall $3'$-subgroup $J$ of $C$ containing $S$. Then $S\cong Q_{16}$ is a Sylow $2$-subgroup of $J$, $3$ does not divide $|J|$, and both $5,7$ divide $|J|$. Every nontrivial $2$-subgroup $U$ of $J$ is cyclic or generalized quaternion. If $U$ is cyclic or $|U|>8$, then $\Aut(U)$ is a $2$-group. For $U\cong Q_8$, we have $\Aut(U)\cong S_4$, so $3\nmid |J|$ again implies that $N_J(U)/C_J(U)$ is a $2$-group. Frobenius's normal $p$-complement theorem~\cite[10.3.2]{Robinson} therefore gives a normal odd-order complement $N$ in $J$.

The group $N$ is a $\{5,7\}$-group containing both primes and no element of order $35$. The proof of Lemma~\ref{orientations} applies: $N$ is Frobenius with a $5$-kernel or a $7$-kernel.

If the kernel is a $5$-group, a Sylow $7$-subgroup of $N$ is cyclic of order $7$. There is an $S$-invariant such subgroup $Q_0$. The map $S\to\Aut(Q_0)\cong C_6$ has abelian image, so it is trivial on $S'\cong C_4$. An element of order $4$ in $S'$ therefore centralizes $Q_0$, yielding order $28$.

If the kernel is a $7$-group $Q_0$, then $Q_0$ is normal in $J$. The involution $z$ is central in $J\le C$, so it centralizes $Q_0$. Choose $T\cong Q_8$ in $S$ containing $z$. Its action on the nonzero elementary abelian $7$-module $M=\Omega_1(Z(Q_0))$ factors through $T/\langle z\rangle\cong C_2^2$.

Every weight of $C_2^2$ on $M$ has a nontrivial kernel. Hence some nonidentity element of $T/\langle z\rangle$ has a nonzero fixed point. A preimage $t\in T$ has order $4$ and centralizes a nonzero vector of $M$, again yielding order $28$.
\end{proof}

\begin{lemma}
\label{vfour}
The group $G$ contains $C_2^2$.
\end{lemma}
\begin{proof}
Let $S\in\Syl_2(G)$. If $S$ contains distinct involutions $u,v$, then $D=\langle u,v\rangle$ is a dihedral $2$-group. If $uv$ has order $2$, the elements $u,v$ already generate $C_2^2$. If $|uv|=2^a\ge4$, the central involution $(uv)^{2^{a-1}}$ commutes with $u$ and is distinct from $u$, again giving $C_2^2$.

Thus the absence of $C_2^2$ implies that $S$ has a unique involution. By the standard description of $p$-groups with a unique subgroup of order $p$~\cite[5.3.6]{Robinson}, $S$ is cyclic or generalized quaternion. The former is excluded by Lemma~\ref{nocyclictwo}, and the latter by Lemma~\ref{noquaterniontwo}. Hence $C_2^2\le S\le G$.
\end{proof}

\subsection{\texorpdfstring{The normal $7$-subgroup and the radical}{The normal 7-subgroup and the radical}}
Throughout this subsection, assume orientation (B), and put
\[
R=O_{7'}(G),\qquad B=G/R,\qquad F=O_7(B),\qquad A\in\Syl_5(G),\quad |A|=5.
\]

\begin{lemma}
\label{sevenfactors}
The group $A$ acts fixed-point-freely on every $7$-chief factor of $G$. If $J=T\rtimes A$ is a Hall $\{3,5\}$-subgroup, then $[T,A]$ centralizes every such factor and $[T,A]\le R$.
\end{lemma}
\begin{proof}
\emph{The position of $7$-factors.} Let $N\trianglelefteq G$ and $5\mid|N|$. Since a full Sylow $5$-subgroup has order $5$, we may conjugate and assume $A\le N$. In the Hall subgroup $Q\rtimes A$, we have $C_Q(A)=1$. Coprime action~\cite[8.2.7(a)]{KurzweilStellmacher} gives $[Q,A]=Q$, and normality gives $Q\le N$. Thus any normal subgroup whose order is divisible by $5$ contains the full $7$-part of $|G|$.

Let $L/K$ be a $7$-chief factor of $G$. If $5\mid|K|$, then $K$ would contain a full Sylow $7$-subgroup and $7$ would not divide $|L:K|$, a contradiction. Hence $K$ is a $5'$-group. If $C_{L/K}(A)\ne0$, Lemma~\ref{localcentralizer} lifts these fixed points to $C_L(A)$. Then $7\mid|C_L(A)|$, so Cauchy's theorem gives an element of order $7$ in $C_L(A)$. It commutes with $A$ and yields order $35$. Therefore $C_{L/K}(A)=0$.

\emph{The Hall $\{3,5\}$-subgroup.} Let $J$ be such a subgroup. For $A\in\Syl_5(J)$, the quotient $N_J(A)/C_J(A)$ embeds in $\Aut(C_5)\cong C_4$, but its order has no prime divisors other than $3,5$. Thus it is trivial. Since $A$ is cyclic, $A\le Z(N_J(A))$. Burnside's normal $p$-complement criterion~\cite[10.1.8]{Robinson} gives a normal $5$-complement $T\in\Syl_3(J)$, and $J=T\rtimes A$.

\emph{Application of Flavell's theorem.} On a $7$-chief factor $M=L/K$, the group $T\rtimes A$ has order coprime to $7$. Its action over a splitting field is completely reducible by Maschke's theorem~\cite[Theorem 1.9]{Isaacs}. Consider each absolutely irreducible summand $W$ and factor out the kernel of the action. The image of $A$ still has order $5$, since it acts fixed-point-freely on $M$ and hence on each nonzero summand. The image of $T$ is a normal $3$-group, and $W$ is a faithful irreducible module for the semidirect product of these images. All hypotheses of Proposition~\ref{extFlavell} now hold for $r=5$ in characteristic $7$. The commutator of the images of $T$ and $A$ is trivial or a nonabelian special $2$-group. As it lies in a $3$-group, the latter is impossible. This holds on each summand, so $[T,A]$ centralizes $M$.

\emph{Passage to $B$.} Put $D=[T,A]$, a $3$-subgroup because $T\trianglelefteq J$. Its image $\bar D$ in $B$ centralizes every $7$-chief factor. Choose a chief series through $R$, the full preimage of $\Phi(F)$, and the full preimage of $F$. The induced series on $F/\Phi(F)$ has only $7$-factors. Matrices acting trivially on every factor of this series are block-unipotent and form a $7$-group. The image of the $3$-group $\bar D$ in it is trivial.

Hence $\bar D$ centralizes $F/\Phi(F)$. By Burnside's theorem on the Frattini quotient~\cite[5.3.2--5.3.3]{Robinson}, the kernel of $\Aut(F)\to\Aut(F/\Phi(F))$ is a $7$-group, so $\bar D$ centralizes $F$ itself. Standard Fitting subgroup properties~\cite[5.4.4]{Robinson} give $\bar D\le C_B(F)\le F$. But $\bar D$ is a $3$-group and $F$ a $7$-group, whence $\bar D=1$. Thus $[T,A]\le R$.
\end{proof}

\begin{lemma}
\label{ninedrop}
We have $9\notin\om(B)$. If $x\in G$ has order $9$, then $\bar x=xR$ has order $3$ and
\[
C_F(\bar x)=1.
\]
\end{lemma}
\begin{proof}
By Lemma~\ref{Hallsections}, the image in $B$ of any Hall $\{3,5\}$-subgroup $J=T\rtimes A$ is Hall. Lemma~\ref{sevenfactors} shows that the images of $T$ and $A$ commute. The group $A$ retains order $5$, since $R$ is a $\{2,3\}$-group. If $B$ had an element of order $9$, it would lie in some such Hall subgroup and commute with its $5$-part, yielding order $45$ in a quotient of $G$. Hence $9\notin\om(B)$.

Let $|x|=9$. By Lemma~\ref{noradicalnine}, $x\notin R$. The order of $\bar x$ divides $9$ and is neither $1$ nor $9$, so it is $3$.

First, $C_F(\bar x)$ cannot contain $E_0\cong C_7^4$. Otherwise $Y=E_0\times\langle\bar x\rangle$ would be a subgroup of $B$. In its full preimage, Hall's theorem~\cite[9.1.7]{Robinson} gives a Hall $\{3,7\}$-subgroup $J$ containing $x$. It maps onto $Y$, and $J\cap R$ is a normal $3$-subgroup. Thus the preimage of the factor $\langle\bar x\rangle$ in $J$ is a normal $3$-subgroup $T$ containing $x$. Schur--Zassenhaus~\cite[9.1.2]{Robinson} gives a complement $E\cong E_0$ to $T$ in $J$.

The group $E$ acts coprimely on $T$. Since $63\notin\Omega$, we have $\exp C_T(e)\le3$ for $e\in E^\#$. Lemma~\ref{rankfour} gives $\exp T\le3$, contrary to $x\in T$.

Now choose a Hall $\{3,5\}$-subgroup $T_0\rtimes A$ containing $x$. By Lemma~\ref{sevenfactors}, $\bar x$ commutes with the image of $A$. Thus $C=C_F(\bar x)$ is $A$-invariant. If $C\ne1$, then $M=\Omega_1(Z(C))$ is a nonzero elementary abelian $7$-module for $A$. The action is fixed-point-free, since a fixed element of order $7$ together with $A$ would yield $35$ in $B$.

Nontrivial irreducible $\F_7C_5$-modules have dimension $\operatorname{ord}_5(7)=4$: $7\equiv2\pmod5$, and $2$ has order $4$ modulo $5$. Hence $\dim_{\F_7}M$ is a positive multiple of $4$, yielding the forbidden subgroup $C_7^4\le C$. Therefore $C_F(\bar x)=1$.
\end{proof}

\begin{lemma}
\label{liftseven}
There is a $7$-subgroup $V\le G$ mapping isomorphically onto $F$ such that
\begin{equation}\label{RV}
[R,V]\le O_2(R).
\end{equation}
The rank of $\Omega_1(Z(V))$ is a positive multiple of $4$.
\end{lemma}
\begin{proof}
\emph{Choosing a lift of $F$.} Apply Lemma~\ref{normalizerlift} with $p=7$ to $R\trianglelefteq G$ and $F\trianglelefteq B$. This gives a complement $V$ in the full preimage of $F$, an isomorphism $V\to F$, and a surjection
\[
N_G(V)\longrightarrow B,\qquad
\ker=C_R(V)=S.
\]
In particular, $S$ is a $\{2,3\}$-group. Choose a subgroup $A\le N_G(V)$ of order $5$; it exists since the normalizer maps onto $B$ and $5\nmid|S|$.

\emph{A compatible element of order $3$.} Take $x\in G$ of order $9$. By Lemma~\ref{ninedrop}, its image $\bar x$ has order $3$, acts fixed-point-freely on $F$, and commutes with a Sylow $5$-subgroup of $B$, namely the image of the $5$-part of a Hall subgroup containing $x$. All Sylow $5$-subgroups of $B$ are conjugate. Conjugating $\bar x$ in $B$ therefore gives $\bar a$ of order $3$ centralizing the image of the chosen $A$, with $C_F(\bar a)=1$.

Apply Lemma~\ref{localcentralizer} to the action of $A$ on $N_G(V)$ and the kernel $S$. It follows that $C_{N_G(V)}(A)$ maps onto $C_B(\bar A)$. In the preimage of $\langle\bar a\rangle$ inside this centralizer, choose a $3$-element $a$ mapping to a generator. Its order cannot be $9$, for together with the centralized subgroup $A$ it would yield $45$. Hence $|a|=3$.

\emph{Action on chief $3$-factors of the radical.} The isomorphism $V\to F$ is compatible with the normalizer action, so $C_V(a)=1$. Also $C_V(A)=1$, since otherwise order $35$ would occur. Put $c=aa_5$, where $a_5$ generates $A$. The elements $a,a_5$ commute, $|c|=15$, and the image of $c$ in $B$ also has order $15$.

Suppose that $V$ acts nontrivially on a chief $3$-factor $M=L/K$ below $R$. Since $7\nmid|M|$, restriction to $V$ is completely reducible and has a nonprincipal constituent. By Lemma~\ref{characters}, all such types have $c$-orbits of length $15$. On the sum of the corresponding isotypic components, Lemma~\ref{norm} gives a nonzero operator $1+c+\cdots+c^{14}$.

The group-theoretic part of Lemma~\ref{norm} applies, since $c$ already retains order $15$ modulo $R$, and hence modulo $L\le R$. It yields order $45$, a contradiction. Thus $V$ centralizes all chief $3$-factors below $R$, and Lemma~\ref{commutator} gives $[R,V]\le O_2(R)$.

Finally, $\Omega_1(Z(V))$ is a nonzero elementary abelian $7$-module on which $A$ acts fixed-point-freely. As in Lemma~\ref{ninedrop}, its dimension is divisible by $\operatorname{ord}_5(7)=4$. Thus $\rk\Omega_1(Z(V))$ is a positive multiple of $4$.
\end{proof}

\begin{lemma}
\label{dichotomy}
We have $R\ne1$ and
\[
O_2(R)=1\quad\text{or}\quad O_3(R)=1.
\]
\end{lemma}
\begin{proof}
Take $x\in G$ of order $9$. Its image has order $3$ by Lemma~\ref{ninedrop}, so $1\ne x^3\in R$ and $R\ne1$.

Suppose that both $O_2(R)$ and $O_3(R)$ are nontrivial. These normal subgroups commute, since their mutual commutators lie in $O_2(R)\cap O_3(R)=1$. The subgroups
\[
M_i=\Omega_1(Z(O_i(R)))\quad(i=2,3)
\]
are nonzero, characteristic in $R$, and hence $V$-invariant. By Lemma~\ref{liftseven}, choose $E\cong C_7^3$ in $V$.

Extend scalars on $M_2$ and $M_3$ to splitting fields for $E$, and choose weights $\lambda_2,\lambda_3$ occurring in the respective modules. Since $E$ has dimension $3$ over $\F_7$ and each weight kernel has codimension at most $1$, their intersection has dimension at least $1$. Choose $e\ne1$ in it.

Thus $C_{M_2}(e)$ and $C_{M_3}(e)$ are nonzero. Take elements $u,v$ of orders $2,3$ in them. The elements $u,v,e$ commute pairwise and have coprime orders $2,3,7$, whence $|uve|=42$, a contradiction.
\end{proof}

\begin{lemma}
\label{noOthree}
The case $O_3(R)=1$ is impossible.
\end{lemma}
\begin{proof}
Suppose that $O_3(R)=1$ and put $N=O_2(R)$. Since $R\ne1$ is solvable, standard Fitting subgroup properties~\cite[5.4.4]{Robinson} give $F(R)=N\ne1$ and $C_R(N)\le N$. A Sylow $2$-subgroup $S$ of $R$ containing $N$ is elementary abelian by Lemma~\ref{twoelementary}. It centralizes $N$, so $S\le C_R(N)\le N$ and $S=N$. Hence $N$ is a full Sylow $2$-subgroup of $R$.

The quotient $U=R/N$ is a nontrivial $3$-group: $x^3\in R$, of order $3$, cannot lie in $N$. Schur--Zassenhaus~\cite[9.1.2]{Robinson} splits $R$ over $N$, with a Sylow $3$-subgroup as complement. Lemma~\ref{noradicalnine} therefore gives $\exp U=3$.

By~\eqref{RV}, $RV/N=(R/N)\times(VN/N)$. The subgroup $VN/N$ is its unique Sylow $7$-subgroup and is therefore normal in $G/N$, since $RV\trianglelefteq G$. It is isomorphic to $V$, since $V\cap N=1$. Choose $E\cong C_7^4$ in its center. The subgroups $U,E$ commute and intersect trivially. The full preimage of $U\times E$ in $G$ has a normal Hall $2$-subgroup $N$. By Schur--Zassenhaus~\cite[9.1.2]{Robinson}, it has an odd-order complement isomorphic to $U\times E$. Choose commuting preimages of $U,E$ inside this complement, and use the same notation for them.

For $u\in U^\#$, we claim that $C_N(u)=1$. Otherwise this centralizer is a nonzero elementary abelian $2$-module invariant under $E$, since $[u,E]=1$. Every weight of $E$ on it has nontrivial kernel: $\rk E=4$, whereas the image of a weight is cyclic. Thus some $e\in E^\#$ fixes a nonzero vector of $C_N(u)$. We obtain three pairwise commuting elements of orders $2,3,7$, hence order $42$.

If $|U|>3$, then $\exp U=3$ implies $C_3^2\le U$: take a central element $z$ of order $3$ and $v\notin\langle z\rangle$. They commute and generate $C_3^2$. Every weight of $C_3^2$ on the nonzero $2$-module $N$ has nontrivial kernel, giving a fixed point of a nonidentity element of $U$, a contradiction. Hence $U\cong C_3$ and $|R|_3=3$.

Take a Hall $\{3,5\}$-subgroup $T\rtimes A$ containing the original $x$ of order $9$. By Lemma~\ref{sevenfactors}, $[T,A]\le R\cap T$. Lemma~\ref{Hallsections} and $|R|_3=3$ give $|R\cap T|=3$. The group $A$ centralizes this intersection, since $|\Aut(C_3)|=2$. Hence $[T,A,A]=1$, while coprime action gives $[T,A,A]=[T,A]$. Thus $[T,A]=1$, and $x$ centralizes $A$, yielding order $45$.
\end{proof}

\begin{lemma}
\label{Rstructure}
If $O_2(R)=1$, then $V\trianglelefteq G$, $[R,V]=1$, and $C_G(V)=RZ(V)\le RV$. Moreover, $R$ has one of the following forms:
\begin{enumerate}
\item[(a)] $R=U$ is a nontrivial $3$-group of exponent $3$;
\item[(b)] $R=U\rtimes\langle z\rangle$, where $U$ is a nonzero elementary abelian $3$-group, $|z|=2$, and $z$ inverts $U$.
\end{enumerate}
\end{lemma}
\begin{proof}
Equation~\eqref{RV} and $O_2(R)=1$ give $[R,V]=1$. The full preimage of $F\trianglelefteq B$ is $RV=R\times V$, normal in $G$. Its Sylow $7$-subgroup $V$ is unique, hence characteristic in $RV$ and normal in $G$.

If $g\in C_G(V)$, its image in $B$ centralizes $F$. Standard Fitting subgroup properties~\cite[5.4.4]{Robinson} place this image in $F$, so $g\in RV$. For $r\in R,v\in V$, the action of $rv$ on $V$ is that of $v$, since $[R,V]=1$. Thus $rv\in C_G(V)$ if and only if $v\in Z(V)$. This gives $C_G(V)=RZ(V)$.

There is no element of order $6$ in $R$, since it would centralize any element of order $7$ in the nontrivial group $V$ and yield $42$. If $R$ has odd order, it is a $3$-group of exponent $3$ by Lemma~\ref{noradicalnine}, giving (a).

Suppose that $R$ has even order. Then $F(R)=O_3(R)\ne1$, since $O_2(R)=1$ and $R$ is solvable. Any involution $z\in R$ acts fixed-point-freely on $F(R)$, since otherwise its centralizer would contain an element of order $3$, giving $6$. By Lemma~\ref{inversionlemma}, $F(R)$ is abelian and $z$ acts by inversion.

A Sylow $2$-subgroup $S$ of $R$ is elementary abelian by Lemma~\ref{twoelementary}. Its action on $F(R)$ is faithful, since $C_R(F(R))\le F(R)$ and $S\cap F(R)=1$. Every nonidentity element induces the same inversion; two distinct involutions would give a nontrivial kernel element. Hence $|S|=2$.

Since $\Aut(C_2)=1$, we have $N_R(S)=C_R(S)$ and $S\le Z(N_R(S))$. Burnside's normal $p$-complement criterion~\cite[10.1.8]{Robinson} gives a normal odd-order complement $U$. It is a full Sylow $3$-subgroup of $R$ and has exponent $3$. The absence of order $6$ makes $z$ fixed-point-free on $U$. Lemma~\ref{inversionlemma} shows that $U$ is abelian and is inverted by $z$, giving (b).
\end{proof}

\begin{lemma}
\label{oddR}
Under the assumptions of Lemma~\ref{Rstructure}(a), $G$ contains no $C_2^2$.
\end{lemma}
\begin{proof}
Here $R=U$ is a nontrivial $3$-group, $V\trianglelefteq G$, $[U,V]=1$, and $C_G(V)=UZ(V)\le UV$. Suppose that $W\cong C_2^2$ lies in $G$.

For any $w\in W^\#$, the conditions $C_U(w)\ne1$ and $C_V(w)\ne1$ cannot both hold. Otherwise choose elements of orders $3,7$ in these centralizers. They commute by $[U,V]=1$, and commute with $w$, yielding $42$. Hence $w$ is fixed-point-free on at least one of $U,V$. By Lemma~\ref{inversionlemma}, such an involutory action is inversion.

At most one involution of $W$ can be fixed-point-free on $V$. If there were two, both would induce inversion and their product, the third nonidentity involution, would centralize $V$. This is impossible because $C_G(V)\le UV$ has odd order. Hence at least two involutions of $W$ invert $U$. Their product $z$ centralizes $U$ and must be fixed-point-free on $V$. Lemma~\ref{inversionlemma} again shows that $V$ is abelian and inverted by $z$. Also $U$ is abelian, being inverted by either of the first two involutions. Hence $C_G(V)=UV$.

Consider any involution $t\in C_G(U)$. It too must be fixed-point-free on $V$, or order $42$ would occur, and hence it inverts $V$. Thus $tz$ centralizes $V$ and lies in $UV$. Write $t=uvz$. The groups $U,V$ commute; $z$ acts trivially on $U$ and inverts $V$. Therefore $t^2=u^2$. Since $t^2=1$ and $U$ is a $3$-group, $u=1$. All involutions centralizing $U$ consequently have the same image $zV$ in $G/V$.

Since $U=R$ is normal in $G$, the property of an involution centralizing $U$ is preserved by conjugation. For any $g\in G$, the involution $z^g$ centralizes $U$, so $z^gV=zV$. Hence $zV\in Z(G/V)$ and $H=V\langle z\rangle\trianglelefteq G$.

Take $x\in G$ of order $9$. It acts on the normal subgroup $H$, and $(|x|,|H|)=(9,2|V|)=1$. Coprime action~\cite[8.2.3(a)]{KurzweilStellmacher} gives an $x$-invariant Sylow $2$-subgroup $\langle t_0\rangle$ of $H$, of order $2$. Its automorphism group is trivial, so $[x,t_0]=1$ and $|xt_0|=18$, a contradiction.
\end{proof}

\begin{lemma}
\label{evenR}
Under the assumptions of Lemma~\ref{Rstructure}(b), a Sylow $2$-subgroup of $G$ is cyclic or generalized quaternion.
\end{lemma}
\begin{proof}
Let $R=U\rtimes\langle z\rangle$ be as in Lemma~\ref{Rstructure}(b), and set $H=C_G(z)$. The group $U=O_3(R)$ is characteristic in $R$ and normal in $G$. The normal subgroup $R/U$ of $G/U$ has order $2$ and is therefore central. Thus $C_{G/U}(zU)=G/U$.

Apply Lemma~\ref{localcentralizer} to the action of $\langle z\rangle$ on $G$ and the normal odd-order subgroup $U$. It gives $G=UH$. Moreover, $U\cap H=C_U(z)=1$, since $z$ inverts the elementary abelian $3$-group $U$. Thus $G=U\rtimes H$ and $z\in Z(H)$. Since $[R,V]=1$, we have $V\le H$. More precisely,
\[
C_H(V)=Z(V)\times\langle z\rangle:
\]
indeed, $C_G(V)=RZ(V)$ and $R\cap H=\langle z\rangle$.

Suppose that $H$ contains an involution $w\ne z$. If $C_V(w)\ne1$, then $C_U(w)=1$, since otherwise order $42$ would occur. By Lemma~\ref{inversionlemma}, $w$ inverts $U$. But $z$ also inverts $U$ and commutes with $w$. Thus the nonidentity involution $wz$ centralizes $U$, and its action on $V$ is that of $w$, since $z$ centralizes $V$. This again yields order $42$. Hence every involution $w\ne z$ in $H$ is fixed-point-free on $V$. If such a $w$ exists, then $V$ is abelian, $w$ acts by inversion, and $C_H(V)=V\times\langle z\rangle$.

Choose $x\in G$ of order $9$ and write it uniquely as $x=ua$, where $u\in U$, $a\in H$. The group $H$ is isomorphic to $G/U$. There is no element of order $9$ in $G/U$: the map to $B=G/R$ has kernel of order $2$, so such an element would retain order $9$, contrary to Lemma~\ref{ninedrop}. Thus $|a|=1$ or $3$. The first is impossible because $\exp U=3$, so $|a|=3$.

On the additive group $U$, put $N_a=1+a+a^2$. The cube formula in $U\rtimes\langle a\rangle$ gives $(va)^3=N_a(v)$ for every $v\in U$. Since $x=ua$ has order $9$, $N_a(u)\ne0$, so $N_a$ is nonzero.

We claim that $M=V\langle z,w\rangle$ is $a$-invariant. Both $V$ and $\langle z\rangle$ are $a$-invariant. The element $w$ acts on the abelian group $V$ as scalar inversion, so $w^a$ induces the same inversion. Thus $(w^a)w\in C_H(V)=V\times\langle z\rangle$, and $w^a\in M$, proving invariance.

The group $M$ has order $4|V|$ and a Sylow $2$-subgroup $\langle z,w\rangle\cong C_2^2$. The action of $\langle a\rangle$ on $M$ is coprime, so~\cite[8.2.3(a)]{KurzweilStellmacher} gives an $a$-invariant $W'\in\Syl_2(M)$. All Sylow $2$-subgroups of $M$ are isomorphic to $C_2^2$. The central involution $z$ lies in every one of them, since $\langle z\rangle\trianglelefteq M$ is a $2$-subgroup. Also $a$ fixes $z$. The stabilizer of a nonzero vector in $\Aut(C_2^2)\cong S_3$ has order $2$, so an automorphism of order $3$ fixing that vector acts trivially on $W'$. Hence $[a,W']=1$.

Choose $w'\in W'\setminus\langle z\rangle$. In characteristic $3$, decompose $U=U_+\oplus U_-$ into the $1$- and $-1$-eigenspaces of $w'$. Since $[a,w']=1$, both summands are invariant under $a$ and $N_a$. The nonzero operator $N_a$ is nonzero on at least one summand. The element $w'$ is trivial on $U_+$, and $w'z$ is trivial on $U_-$, since $z$ acts as $-1$ on all of $U$.

Choose $v$ in the appropriate summand with $N_a(v)\ne0$. The element $va$ has order $9$, since its cube is a nonzero element of $U$ of order $3$. It centralizes the corresponding involution $w'$ or $w'z$, which fixes $v$ and commutes with $a$. This yields order $18$.

Thus $H$ has no involution other than $z$. The standard description of $p$-groups with a unique subgroup of order $p$~\cite[5.3.6]{Robinson} shows that a Sylow $2$-subgroup of $H$ is cyclic or generalized quaternion. Since $|G:H|=|U|$ is odd, it is also a Sylow $2$-subgroup of $G$.
\end{proof}

\begin{lemma}
\label{sevenconclusion}
If $G$ contains $C_2^2$, orientation \textup{(B)} is impossible.
\end{lemma}
\begin{proof}
By Lemma~\ref{dichotomy}, either $O_2(R)=1$ or $O_3(R)=1$. If $O_3(R)=1$, Lemma~\ref{noOthree} gives a contradiction regardless of $O_2(R)$. Thus $O_2(R)=1$ and $O_3(R)\ne1$.

Apply Lemma~\ref{Rstructure}. If $R$ has odd order, Lemma~\ref{oddR} excludes it using the assumed subgroup $C_2^2$. If $R$ has even order, it has the form $U\rtimes C_2$ described there, and Lemma~\ref{evenR} shows that a Sylow $2$-subgroup of $G$ is cyclic or generalized quaternion. Both have a unique involution~\cite[5.3.6]{Robinson} and contain no $C_2^2$, again contradicting the hypothesis. These exhaust the possibilities.
\end{proof}

\section{\texorpdfstring{The branch with a $5$-kernel: successive reductions}{The branch with a 5-kernel: successive reductions}}
Assume orientation (A) throughout this section, and put
\[
R=O_{5'}(G),\qquad B=G/R,\qquad P=O_5(B).
\]
We establish the following structural reduction:
\begin{equation}\label{Asetup}
P\in\Syl_5(B),\quad B=P\rtimes D,\quad C_D(P)=1,
\quad Q\cong C_7\trianglelefteq D,\quad C_P(Q)=1.
\end{equation}
Here $D$ complements $P$ in $B$, and $Q$ is its Sylow $7$-subgroup. 

\begin{lemma}
\label{faithfullift}
Let $L$ be a finite solvable group with $O_p(L)=1$, and let $M$ be a faithful module in characteristic $p$. Then the semisimple module formed from the composition factors of $M$ is faithful and has a faithful complex lift. For a $p$-regular element $q\in L$, the lift preserves the multiplicity of the eigenvalue $1$.
\end{lemma}
\begin{proof}
Extend scalars to an algebraic closure. The kernel of the action on the direct sum of the composition factors of $M$ consists of transformations acting trivially on each factor of a composition series. In a basis adapted to that series they are unipotent. Since the original action on $M$ is faithful, this kernel is a normal $p$-subgroup of $L$ and hence is trivial.

By the Fong--Swan theorem~\cite[Corollary 0.33]{ManzWolf}, every simple module for a solvable group lifts to characteristic zero. The solvable group $L$ is $p$-solvable for every prime $p$, so the hypothesis of that theorem is satisfied. Take the direct sum of the corresponding complex representations, and denote its character by $\chi$. Its restriction to the $p$-regular elements is the Brauer character of the semisimple module.

We verify faithfulness separately. If a $p$-regular element $h$ is in the kernel of the lift, then $\chi(h)=\chi(1)$. Thus its Brauer character value equals the module dimension. This value is a sum of complex roots of unity obtained by lifting the eigenvalues of $h$. Equality of the sum to the number of terms forces every root to equal $1$. Since $h$ acts semisimply, it acts trivially on the semisimple module. Faithfulness of that module gives $h=1$.

Hence the kernel of the complex lift contains no nonidentity $p'$-elements. By Cauchy's theorem it is a $p$-group, and normality together with $O_p(L)=1$ makes it trivial.

Finally, on the cyclic $p'$-subgroup $\langle q\rangle$, the Brauer character is an ordinary character after identifying roots of unity. Its values on all powers of $q$ determine eigenvalue multiplicities, so the multiplicity of $1$ is unchanged by lifting. These multiplicities are also unchanged when passing from $M$ to its composition factors, since $\langle q\rangle$ acts semisimply.
\end{proof}

\begin{lemma}
\label{shultseven}
Let a finite solvable group $L$ have a faithful complex representation $M$. If $q\in L$ has order $7$ and $C_M(q)=0$, then $q\in O_7(L)$.
\end{lemma}
\begin{proof}
Apply Proposition~\ref{extEspuelas} with $X=L$, $K=\mathbb C$, $p=7$, $n=1$, and $x=q$. Solvability implies $7$-solvability. The representation is faithful, characteristic zero does not divide $|O_{7'}(L)|$, and $C_M(q)=0$ by hypothesis.

The prime $7$ is not Fermat: $7-1=6$ is not a power of $2$. Thus the extra hypothesis on abelian Sylow $2$-subgroups is unnecessary, and the condition for $p=2$ does not arise. The conclusion of Proposition~\ref{extEspuelas} is $q^{7^{1-1}}=q\in O_7(L)$.
\end{proof}

\begin{lemma}
\label{Asetupderived}
Conditions~\eqref{Asetup} hold in orientation \textup{(A)}.
\end{lemma}
\begin{proof}
By Lemma~\ref{radical}, $F(B)=P\ne1$ and $C_B(P)\le P$. Consider the action of $B$ on $M=P/\Phi(P)$. Its kernel $J_0$ contains $P$, since inner automorphisms of $P$ act trivially on $M$. The image of $J_0$ in $\Aut(P)$ lies in the kernel of $\Aut(P)\to\Aut(P/\Phi(P))$, a $5$-group. The kernel of the action of $J_0$ on $P$ lies in $C_B(P)\le P$. Hence $J_0$ is a normal $5$-subgroup of $B$, and $J_0=P$.

Thus $L=B/P$ acts faithfully on $M$. Also $O_5(L)=1$: the preimage of a nontrivial normal $5$-subgroup of $L$ would be a normal $5$-subgroup of $B$ larger than $P$.

In (A), a full Sylow $7$-subgroup of $G$ has order $7$. Since $R$ is a $\{2,3\}$-group and $P$ a $5$-group, Sylow $7$-subgroups of both $B$ and $L$ also have order $7$. Choose $q\in B$ of order $7$. The absence of $35$ in $B$ gives $C_P(q)=1$, hence $C_M(q)=0$.

Apply Lemma~\ref{faithfullift} to $L$ with $p=5$. This gives a faithful complex representation in which the image of $q$ remains fixed-point-free. By Lemma~\ref{shultseven}, this image lies in $O_7(L)$. Thus the Sylow $7$-subgroup $Q_L\cong C_7$ is normal in $L$.

We claim that $5\nmid |L|$. Otherwise $L$ contains $a$ of order $5$. It normalizes $Q_L$ and acts trivially on it, since $5\nmid|\Aut(C_7)|=6$. Then $L$ would contain order $35$, impossible in a quotient of $G$. Hence $P$ is a full Sylow $5$-subgroup of $B$.

Schur--Zassenhaus gives $B=P\rtimes D$, where $D\cong L$. The preimage of $Q_L$ in $D$ is a normal $Q\cong C_7$. We have $C_D(P)\le D\cap P=1$ since $C_B(P)\le P$, and $C_P(Q)=1$ since order $35$ is absent. All conditions in~\eqref{Asetup} follow.
\end{proof}

\begin{lemma}
\label{noeightD}
Under~\eqref{Asetup}, $D$ has no element of order $8$.
\end{lemma}
\begin{proof}
Suppose that $y\in D$ has order $8$. Normality of $Q$ gives a homomorphism $\langle y\rangle\to\Aut(Q)\cong C_6$. Its image is a $2$-subgroup of $C_6$, of order $1$ or $2$. Hence $y^2$ centralizes $Q$. But $|y^2|=4$ and $|Q|=7$, so $y^2q$ has order $28$ for $q\in Q^\#$, contrary to $\om(D)\subseteq\Omega$.
\end{proof}

\begin{lemma}
\label{cubeA}
Under~\eqref{Asetup}, if $x\in D$ has order $9$, then $C_P(x^3)=1$.
\end{lemma}
\begin{proof}
The homomorphism $\langle x\rangle\to\Aut(Q)\cong C_6$ has image of order $1$ or $3$. In the first case $x$ would commute with a generator of $Q$, yielding $63$. Thus the image has order $3$, and $x^3$ centralizes $Q$.

Suppose that $C=C_P(x^3)\ne1$. Both $Q$ and $x$ normalize $P$ and centralize $x^3$, so they normalize $C$. The group $Q$ is fixed-point-free on the nonzero $5$-module $M=C/\Phi(C)$. Indeed, $C_C(Q)=1$ follows from $C_P(Q)=1$, and coprime action~\cite[8.2.2(a)]{KurzweilStellmacher} gives $C_M(Q)=C_C(Q)\Phi(C)/\Phi(C)=0$.

The element $x^3$ acts trivially on $C$, hence on $M$, so the operator $x$ on $M$ satisfies $x^3=1$. Over a splitting field all weights of $Q$ are nontrivial. Conjugation by $x$ on $Q$ has order $3$ and permutes the six nontrivial weights in cycles of length $3$. Every occurring weight belongs to such a cycle, whose weight spaces have an $x$-invariant sum. Lemma~\ref{norm} gives a nonzero fixed vector of $x$ on that sum. Hence $C_M(x)\ne0$.

Since $(9,|C|)=1$, coprime action~\cite[8.2.2(a)]{KurzweilStellmacher} lifts this fixed point to $C_C(x)\ne1$. Choose $u$ of order $5$ in that centralizer. Then $ux=xu$ and $|ux|=45$, a contradiction.
\end{proof}

\begin{lemma}
\label{normalfive}
Under~\eqref{Asetup}, suppose that $P$ lifts to a normal Sylow $5$-subgroup $\widetilde P\trianglelefteq G$, that $R$ centralizes $\widetilde P$, and that $9\in\om(D)$. If the Sylow $2$-subgroups of $D$ have order at most $2$, then $G$ does not satisfy~\eqref{spectrum}.
\end{lemma}
\begin{proof}
Choose $g\in G$ of order $8$. Its image in $B$ is a $2$-element, hence lies in a complement conjugate to $D$: this follows from Hall's theorem~\cite[9.1.7]{Robinson} for the set of primes other than $5$. Conjugating $g$, assume $gR=d\in D$. By hypothesis $|D|_2\le2$, so $|d|\le2$. Since $[R,\widetilde P]=1$, the actions of $g$ on $\widetilde P$ and of $d$ on $P$ correspond under $\widetilde P\to P$.

If $C_{\widetilde P}(g)\ne1$, choose $u$ of order $5$ in it; then $|ug|=40$, a contradiction. Thus $C_{\widetilde P}(g)=1$. In particular, $d\ne1$ and its action has order $2$. Lemma~\ref{inversionlemma} shows that $P$ is abelian and $d$ acts by inversion.

Inversion commutes with every automorphism of an abelian group: $(u^{-1})^\beta=(u^\beta)^{-1}$ for $\beta\in\Aut(P)$. Since the action of $D$ on $P$ is faithful by~\eqref{Asetup}, we have $d\in Z(D)$. By hypothesis, choose $v\in D$ of order $9$. The commuting elements $d,v$ have coprime orders $2,9$, so $|dv|=18$. This is impossible since $D\le B$ and $\om(B)\subseteq\Omega$.
\end{proof}

\begin{lemma}
\label{orbit21}
Assume~\eqref{Asetup}. Suppose that a subgroup $\widetilde P\cong P$ mapping onto $P$ and an element $c$ of order $21$ in $N_G(\widetilde P)$ have been found, such that the parts of $c$ of orders $3,7$ act fixed-point-freely on $\widetilde P$ and its image in $B$ has order $21$. Then $[R,\widetilde P]=1$, and $\widetilde P\trianglelefteq G$ is a full Sylow $5$-subgroup of $G$.
\end{lemma}
\begin{proof}
Write $c=ab=ba$, where $|a|=3$, $|b|=7$. By hypothesis, $C_{\widetilde P}(a)=C_{\widetilde P}(b)=1$. Lemma~\ref{characters} shows that no power of $c$ generating the subgroup of order $3$ or $7$ fixes a nonprincipal irreducible type of $\widetilde P$ in characteristic $2$ or $3$. Thus the stabilizer of such a type in $\langle c\rangle$ is trivial and its orbit has length $21$.

Consider any $G$-chief factor $M=L/K$ below $R$. It is an elementary abelian $r$-group, with $r=2$ or $3$, since $G$ is solvable and $R$ a $\{2,3\}$-group. Suppose that $\widetilde P$ acts nontrivially on $M$. Restriction to $\widetilde P$ is completely reducible by Maschke's theorem~\cite[Theorem 1.9]{Isaacs}, since $r\ne5$. Over a splitting field it has a nonprincipal constituent, and $c$ cyclically permutes the $21$ corresponding isotypic components. Lemma~\ref{norm} gives a nonzero operator $1+c+\cdots+c^{20}$ on $M$.

Apply the group-theoretic part of that lemma to $L/K$ and $c$. The image of $c$ modulo $L$ has order $21$, since it already has order $21$ modulo $R$ and $L\le R$. We obtain an element whose order is divisible by $42$ if $r=2$, and by $63$ if $r=3$. Both contradict the spectrum.

Thus $\widetilde P$ centralizes every chief factor below $R$. Lemma~\ref{commutator} gives $[R,\widetilde P]=1$. The full preimage of the normal subgroup $P\trianglelefteq B$ is $R\widetilde P\trianglelefteq G$. Since $R\cap\widetilde P=1$ and the factors commute, $R\widetilde P=R\times\widetilde P$. Hence $\widetilde P$ is the unique Sylow $5$-subgroup of $R\widetilde P$, is characteristic in it, and is normal in $G$.

Finally, $5\nmid|R|$ and $|B|_5=|P|$, so $|G|_5=|\widetilde P|$.
\end{proof}

\subsection{\texorpdfstring{The case $9\in\omega(B)$}{The case of elements of order 9 in B}}
Suppose that $9\in\om(B)$, and hence $9\in\om(D)$. By Lemma~\ref{cubeA}, the cube of the corresponding element acts fixed-point-freely on $P$ and centralizes $Q$. Before applying Lemma~\ref{orbit21}, we justify the lift of an element of order $21$.

\begin{lemma}
\label{lift21}
Assume~\eqref{Asetup}. If $B$ contains $Q_0\times\langle a_0\rangle\cong C_7\times C_3$ whose two prime-order factors act fixed-point-freely on $P$, then there exist $\widetilde P\cong P$ over $P$ and $c\in N_G(\widetilde P)$ satisfying all hypotheses of Lemma~\ref{orbit21}.
\end{lemma}
\begin{proof}
Choose a complement $\widetilde P$ to $R$ in the preimage of $P$. Conjugacy of complements gives $G=R N_G(\widetilde P)$. Put $N=N_G(\widetilde P)$ and $S=N\cap R$. Then $N/S\cong B$ and $S=C_R(\widetilde P)$: an element of $R$ normalizing $\widetilde P$ centralizes it, since its commutators with $\widetilde P$ lie in $R\cap\widetilde P=1$.

In the preimage of $Q_0\times\langle a_0\rangle$ inside $N$, choose a Hall $\{3,7\}$-subgroup $J$. Its index is a power of $2$, so it maps onto the entire group of order $21$. The group $T=J\cap S$ is a normal $3$-subgroup, and $J/T\cong C_{21}$.

The preimage of $Q_0$ in $J$ has the form $T\rtimes Q_1$, where $|Q_1|=7$, and is normal in $J$. Conjugacy of complements gives $J=T N_J(Q_1)$, so $N_J(Q_1)$ maps onto $J/T$. Since $J/T$ is abelian, $N_J(Q_1)$ centralizes $Q_1$: a commutator with $Q_1$ lies in both $Q_1$ and $T$.

Take a $3$-element $a\in N_J(Q_1)$ whose image generates $\langle a_0\rangle$. Its order is $3$: order $9$ together with the centralized subgroup $Q_1$ would yield $63$, and larger powers of $3$ are already forbidden. Thus $c=qa$, for a generator $q\in Q_1$, has order $21$ and retains this order in $B$. The actions of $q,a$ on $\widetilde P$ correspond to those of their images on $P$, since $\widetilde P\to P$ is an isomorphism. Both are fixed-point-free, as required.
\end{proof}

We first classify the solvable groups whose element orders belong to $\{1,2,3\}$.

\begin{lemma}
\label{Estructure}
Let $E$ be a finite solvable group with $\om(E)\subseteq\{1,2,3\}$. Then one of the following holds:
\begin{enumerate}
\item[(1)] $E$ is a $3$-group of exponent at most $3$;
\item[(2)] $E=U_2$ is a nonzero elementary abelian $2$-group;
\item[(3)] $E=U_2\rtimes\langle t\rangle$, where $U_2$ is a nonzero elementary abelian $2$-group, $|t|=3$, and $C_{U_2}(t)=1$;
\item[(4)] $E=U\rtimes\langle z\rangle$, where $U$ is a nonzero elementary abelian $3$-group, $|z|=2$, and $z$ inverts $U$.
\end{enumerate}
Case \textup{(1)} includes the trivial group.
\end{lemma}
\begin{proof}
Cauchy's theorem gives $\pi(E)\subseteq\{2,3\}$. If $|E|$ is odd, we have (1). Suppose that $|E|$ is even. Every Sylow $2$-subgroup is elementary abelian, since order $4$ is absent.

If both $O_2(E),O_3(E)$ were nontrivial, they would commute and yield order $6$. Thus at most one is nontrivial. The Fitting subgroup of a nontrivial solvable group is nontrivial.

Let $F(E)=O_2(E)=U_2\ne1$. If $S\in\Syl_2(E)$ contains $U_2$, then $S$ is abelian, so $S\le C_E(U_2)\le U_2$. Hence $S=U_2$, and $U_2\trianglelefteq E$ is a full Sylow $2$-subgroup. Schur--Zassenhaus gives $E=U_2\rtimes T$, where $T$ is a $3$-group of exponent at most $3$. Every nonidentity element of $T$ acts fixed-point-freely on $U_2$, or order $6$ would occur. If $|T|>3$, take $c\in Z(T)$ of order $3$ and $v\in T\setminus\langle c\rangle$. Since $\exp T\le3$, they generate $C_3^2$. Over a splitting field, each weight of this subgroup on the nonzero $2$-module $U_2$ has nontrivial kernel, contradicting the fixed-point-free action of every nonidentity element. Hence $|T|\le3$, giving (2) or (3).

Let $F(E)=O_3(E)\ne1$. Any involution $z\in E$ is fixed-point-free on $F(E)$, so Lemma~\ref{inversionlemma} makes $F(E)$ abelian and $z$ acts by inversion. A Sylow $2$-subgroup $S$ acts faithfully on $F(E)$ because $C_E(F(E))\le F(E)$. All its nonidentity elements induce the same inversion, so $|S|=2$.

For $S=\langle z\rangle$ of order $2$, we have $N_E(S)=C_E(S)$ since $\Aut(C_2)=1$. Thus $S\le Z(N_E(S))$, and Burnside's normal $p$-complement criterion~\cite[10.1.8]{Robinson} gives a normal odd-order complement $U$. It is a $3$-group on which $z$ is fixed-point-free. By Lemma~\ref{inversionlemma}, $U$ is abelian and inverted by $z$. Since order $9$ is absent, $U$ is elementary abelian, giving (4).
\end{proof}

\begin{lemma}
\label{doublefree}
Let $U_2$ be a normal elementary abelian $2$-subgroup of $U_2\langle t\rangle$, where $|t|=3$ and $C_{U_2}(t)=1$. Suppose that $U_2\langle t\rangle$ acts on a finite $5$-group $P$ with $C_P(t)=1$. Then $[U_2,P]=1$.
\end{lemma}
\begin{proof}
The conclusion is clear for $P=1$. Suppose that $P\ne1$ and pass to $M=P/\Phi(P)\ne0$. Since $3\nmid|P|$, coprime action~\cite[8.2.2(a)]{KurzweilStellmacher} gives $C_M(t)=0$ from $C_P(t)=1$.

Extend scalars to a splitting field for the elementary abelian $2$-group $U_2$. By Maschke's theorem~\cite[Theorem 1.9]{Isaacs}, $M$ is a direct sum of weight spaces for $U_2$. The action of $t$ on these spaces agrees with its permutation of weights. Since $C_{U_2}(t)=1$, Lemma~\ref{characters} shows that every nontrivial weight has orbit of length $3$. If such a weight occurs in $M$, Lemma~\ref{norm} makes $1+t+t^2$ nonzero on the sum of the three weight spaces. Its image is fixed by $t$, since $t^3=1$, contrary to $C_M(t)=0$.

Thus only the principal weight occurs, and $U_2$ centralizes $M$. By Burnside's theorem on the Frattini quotient~\cite[5.3.2--5.3.3]{Robinson}, the kernel of $\Aut(P)\to\Aut(M)$ is a $5$-group. The image of $U_2$ is both a $2$-group and a subgroup of this kernel, so it is trivial. Hence $[U_2,P]=1$.
\end{proof}

\begin{lemma}
\label{centralizerE}
Under~\eqref{Asetup},
\[
H_0=C_D(Q)=Q\times E,\qquad E\trianglelefteq D,\qquad
\om(E)\subseteq\{1,2,3\},\qquad D/H_0\le C_6.
\]
If $|E|$ is odd, then $|D|_2\le2$. If $|E|$ is even and $9\in\om(D)$, case \textup{(2)} of Lemma~\ref{Estructure} is impossible, while case \textup{(3)} forces $O_2(E)=1$ and is therefore also impossible.
\end{lemma}
\begin{proof}
The group $Q$ is central in $H_0$ and is a Hall $7$-subgroup of $H_0$. A complement $E$ exists and is unique: complements are conjugate by elements of $Q$, but $Q\le Z(H_0)$. Hence $H_0=Q\times E$. The subgroup $Q$ is characteristic in $H_0$ as its unique Sylow $7$-subgroup, and $E$ is characteristic as the unique complement. Since $H_0=C_D(Q)\trianglelefteq D$, we have $E\trianglelefteq D$. The action on $Q$ gives $D/H_0\le\Aut(Q)\cong C_6$. The prime divisors of $D$ lie in $\{2,3,7\}$ and $7\nmid|E|$, so $\pi(E)\subseteq\{2,3\}$. Orders $4,6,9$ in $E$ would yield $28,42,63$ in $H_0$, respectively. Every other order greater than $3$ involving only primes $2,3$ is divisible by $4$, $6$, or $9$, so $\om(E)\subseteq\{1,2,3\}$. If $E$ has odd order, $|D|_2=|D/H_0|_2\le2$.

Let $x\in D$ have order $9$. Its action on $Q$ has order $3$, or order $63$ would occur. Hence $x^3\in H_0$, and in $Q\times E$ this element lies in $E$ and has order $3$. This excludes (2).

In case (3), put $U_2=O_2(E)$. It is normal in $D$, and $t=x^3\in E$ acts fixed-point-freely on it, since $6\notin\om(E)$. Also $C_P(t)=1$ by Lemma~\ref{cubeA}. Lemma~\ref{doublefree} gives $[U_2,P]=1$, and $C_D(P)=1$ then gives $U_2=1$, contrary to (3).
\end{proof}

\begin{lemma}
\label{gapD}
Under~\eqref{Asetup}, if $9\in\om(D)$, then $|D|_2\le2$.
\end{lemma}
\begin{proof}
By Lemmas~\ref{Estructure} and~\ref{centralizerE}, it suffices to consider
\[
E=U\rtimes\langle z\rangle,\qquad
U\cong C_3^d,\quad d\ge1,\quad u^z=u^{-1}.
\]
The group $U=O_3(E)$ is normal in $D$. Choose a Hall $\{2,3\}$-subgroup $J\le D$. It contains $U$, and $E\cap J=U\rtimes\langle z\rangle$ after choosing $z$ in this intersection. We have $J/(E\cap J)\le C_6$, and $(E\cap J)/U\cong C_2$ is central in $J/U$. Thus $J/U$ is a central cyclic extension of a cyclic group and is abelian. Moreover,
\[
|J/U|_3\le3,\qquad |J|_2\le4.
\]
A Sylow $3$-subgroup $T$ of $J$ is normal, contains $U$, and satisfies $|T:U|=3$, since $T$ has an element of order $9$ whereas $\exp U=3$.

Choose $x\in T$ of order $9$, and write $c=x^3\in U^\#$, using additive notation on $U$. Let $\alpha$ be conjugation by $x$ on $U$. Since $x^3\in U$ and $U$ is abelian, $\alpha^3=1$. Set
\[
N=1+\alpha+\alpha^2=(\alpha-1)^2,\qquad W=N(U).
\]
Since $J/U$ is abelian, $x^z=ux$ for some $u\in U$. On the one hand, $(x^z)^3=c^z=-c$; on the other, $(ux)^3=c+N(u)$. In characteristic $3$, this gives $N(u)=-2c=c$. Thus $c\in W$ and $W\ne0$.

Every $w\in W$ has the form $c+N(v)$ and hence equals $(vx)^3$ for some $v\in U$. If $w\ne0$, the element $vx$ has order $9$, so Lemma~\ref{cubeA} gives $C_P(w)=1$. All nonidentity elements of the elementary abelian group $W$ are therefore fixed-point-free on $P$. This is possible only if $\rk W=1$; otherwise each weight of $W$ on $P/\Phi(P)\ne0$ would have nontrivial kernel. Hence
\[
W\cong C_3.
\]

The subgroup $W$ is $J$-invariant. Indeed, $J$ normalizes $T,U$, and its action on the cyclic group $T/U$ sends $xU$ to $x^\varepsilon U$, with $\varepsilon=1$ or $2$. It therefore conjugates $\alpha$ to $\alpha^\varepsilon$. Since $\alpha^3=1$, in characteristic $3$ we have
\[
1+\alpha^2+\alpha^4=1+\alpha+\alpha^2=N,
\]
so the image of the norm operator is preserved.

Suppose that $J$ contains $y$ of order $4$. Its square acts trivially on $Q$ and is an involution of $E$. Thus $y^2$ inverts $U$, and hence $W$. But the square of every automorphism of $W\cong C_3$ is trivial, since $\Aut(C_3)\cong C_2$. This is a contradiction, so $J$ has no element of order $4$.

It remains to exclude $S\cong C_2^2$ in $J$. Such an $S$ has an involution $s\ne1$ centralizing $W$, since $\Aut(W)\cong C_2$. Since $J/U$ is abelian, $s$ centralizes $T/U$. Coprime lifting of centralizers gives $b\in C_T(s)$ with $bU=xU$. The action of $b$ on $U$ is $\alpha$, and $b^3\in W$.

The operators $s$ and $\alpha$ on $U$ commute, since $x^s$ differs from $x$ by an element of the abelian group $U$. Thus $s$ commutes with $N$. Since $s$ is trivial on $W=N(U)$ and $2$ is invertible in $\F_3$, averaging gives
\[
N(C_U(s))=W.
\]
Indeed, for $w=N(u)\in W$, the element $u_0=2^{-1}(u+u^s)$ lies in $C_U(s)$ and satisfies $N(u_0)=w$, since $s$ fixes $W$. Choose $v\in C_U(s)$ with $b^3+N(v)\ne0$, possible because $N(C_U(s))=W\ne0$. Then $(vb)^3=b^3+N(v)$, and $vb$ has order $9$ and centralizes the involution $s$. This yields the forbidden order $18$.

We know $|J|_2\le4$. If $|J|_2=4$, its Sylow $2$-subgroup would be $C_4$ or $C_2^2$, both just excluded. Thus $|J|_2\le2$. Since $J$ is Hall in $D$, $|D|_2=|J|_2\le2$.
\end{proof}

\begin{corollary}[Excluding order $9$ in the quotient]\label{A9survive}
Under~\eqref{Asetup}, we have $9\notin\om(B)$.
\end{corollary}
\begin{proof}
If $x\in D$ has order $9$, its cube centralizes $Q$ and acts fixed-point-freely on $P$ by Lemma~\ref{cubeA}. The action of $Q$ is also fixed-point-free. Lemmas~\ref{lift21} and~\ref{orbit21} give a normal full Sylow $5$-subgroup $\widetilde P$ of $G$ centralized by $R$. By Lemma~\ref{gapD}, $|D|_2\le2$. Lemma~\ref{normalfive} now gives a contradiction.
\end{proof}

\subsection{\texorpdfstring{The drop from order $9$ to order $3$}{The drop from order 9 to order 3}}
By Lemma~\ref{noradicalnine} and Corollary~\ref{A9survive}, every $x\in G$ of order $9$ has image $\bar x$ of order $3$ in $B$. After conjugation assume $\bar x\in D$. Its action on $Q$ is trivial or has order $3$. We first justify compatible lifts in the nontrivial case.

\begin{lemma}
\label{gaplift}
Assume~\eqref{Asetup}. There is no $x\in G$ of order $9$ such that $a_0=xR\in D$ has order $3$ and $[a_0,Q]\ne1$.
\end{lemma}
\begin{proof}
Since $\exp P=5$, its center $Z(P)$ is elementary abelian and is a nonzero $\F_5[Q\rtimes\langle a_0\rangle]$-module. The action of $Q$ on it is fixed-point-free, so $6\mid\dim_{\F_5}Z(P)$. Over a splitting field, $a_0$ has orbits of length $3$ on the nontrivial weights of $Q$. Since $a_0^3=1$, the fixed-point subspace on the sum belonging to each orbit has one third of the dimension of that sum. Hence
\[
\dim_{\F_5}C_{Z(P)}(a_0)=\tfrac13\dim_{\F_5}Z(P)\ge2.
\]
Choose $A_0\cong C_5^2$ in $C_{Z(P)}(a_0)$.

As in Lemma~\ref{lift21}, choose a complement $\widetilde P$ to $R$ in the preimage of $P$, and put $N=N_G(\widetilde P)$. The map $N\to B$ is surjective and has kernel $S=C_R(\widetilde P)$, centralizing $\widetilde P$. Let $A\le Z(\widetilde P)$ be the preimage of $A_0$ under the isomorphism $\widetilde P\to P$.

Take a $3$-element $a\in N$ with image $a_0$. It centralizes $A$, since the actions of $a$ on $\widetilde P$ and of $a_0$ on $P$ correspond. Its order is $3$, since order $9$ together with $A$ would yield $45$.

Suppose that $\widetilde P$ acts nontrivially on a chief $3$-factor $M=L/K$ below $R$. Work over a splitting field for $\widetilde P$, and choose a nonprincipal irreducible type $\theta$ occurring in $M$. A generator of $Q$ lifts in $N$ to $q$ of order $7$. In $B$, we have $\bar q^{a_0}=\bar q^k$ with $k=2$ or $4$ modulo $7$, since $a_0$ induces a nontrivial automorphism of order $3$ on $Q$. Hence $q^a q^{-k}\in S$. The subgroup $S=C_R(\widetilde P)$ acts trivially by conjugation on $\widetilde P$ and on its irreducible types. Thus the permutations of types induced by $a,q$ satisfy $q^a=q^k$.

By Lemma~\ref{characters}, the $q$-orbit of $\theta$ has length $7$. The element $a$ cannot fix every type in this orbit: otherwise applying $a$ to $\theta^q$ would give $\theta^q=\theta^{q^k}$, contrary to the orbit length. Thus some occurring type has an $a$-orbit of length $3$.

Since $A\le Z(\widetilde P)$, it acts on the chosen isotypic component by a scalar character. This character has nontrivial kernel on $A\cong C_5^2$. Choose $p\in A^\#$ in that kernel. The element $a$ centralizes $p$, so $p$ acts trivially on all three components of the chosen $a$-orbit. Lemma~\ref{norm} shows that $1+a+a^2$ is nonzero on $C_M(p)$.

The subgroups $K,L$ are normal in $G$ and lie in the $\{2,3\}$-group $R$, so the action of $p$, of order $5$, on $L$ is coprime. Coprime action~\cite[8.2.2(a)]{KurzweilStellmacher} maps $C_L(p)$ onto $C_M(p)$. The map is compatible with $a$, since $[a,p]=1$. Apply the group-theoretic part of Lemma~\ref{norm}, with $r=m=3$, to the section $C_L(p)/C_K(p)\cong C_M(p)$ and the element $a$. Its image modulo $C_L(p)$ has order $3$, since $a$ retains this order modulo $R$. We obtain $v\in C_L(p)$ such that $9\mid |va|$. Since $va$ centralizes $p$, an order divisible by $45$ occurs, a contradiction.

Thus $\widetilde P$ centralizes all chief $3$-factors below $R$. Lemma~\ref{commutator} gives
\[
[R,\widetilde P]\le O_2(R)=R_2.
\]
In $G/R_2$, the image of $\widetilde P$ is normal: it centralizes $R/R_2$ and is the unique Sylow $5$-subgroup of the normal subgroup $R\widetilde P/R_2$. Conjugation by the original $x$ on this image agrees with the action of $a_0$, since the kernel $R/R_2$ centralizes it. Hence $xR_2$ centralizes the image of $A$. Its order is $9$, since $R_2$ is a $2$-group. This yields order $45$ in $G/R_2$.
\end{proof}

\begin{lemma}
\label{centraldrop}
Assume~\eqref{Asetup}. If $|x|=9$, $a_0=xR\in D$ has order $3$, and $[a_0,Q]=1$, then $C_P(a_0)=1$.
\end{lemma}
\begin{proof}
Suppose that $C=C_P(a_0)\ne1$. Since $[a_0,Q]=1$, the group $Q$ preserves $C$. It acts fixed-point-freely on $C$, because $C_C(Q)\le C_P(Q)=1$. The center $Z(C)$ is nonzero, $Q$-invariant, and has exponent $5$, since $\exp P=5$.

The minimum degree of an extension of $\F_5$ containing a nontrivial seventh root of unity is $\operatorname{ord}_7(5)=6$. Thus every nontrivial simple $\F_5Q$-module has dimension $6$. There are no principal constituents in $Z(C)$, so its dimension is a positive multiple of $6$. Choose $E_0\cong C_5^4$ in $Z(C)\le C_P(a_0)$. Then $E_0\times\langle a_0\rangle$ is a subgroup of $B$.

Let $Y$ be its full preimage in $G$; it contains the original element $x$. By Hall's theorem~\cite[9.1.7]{Robinson}, choose a Hall $\{3,5\}$-subgroup $J$ of $Y$ containing $\langle x\rangle$. Lemma~\ref{Hallsections} shows that $J$ maps onto $Y/R=E_0\times\langle a_0\rangle$. Also $J\cap R$ is a normal $3$-subgroup of $J$, since $R$ is a $\{2,3\}$-group.

Let $T$ be the preimage of the factor $\langle a_0\rangle$ in $J$. Then $T\trianglelefteq J$ is a $3$-group, $x\in T$, and $J/T\cong E_0$. Schur--Zassenhaus~\cite[9.1.2]{Robinson} gives $J=T\rtimes E$ with $E\cong C_5^4$. For every $e\in E^\#$, $C_T(e)$ has no element of order $9$, or together with $e$ it would yield $45$. Larger powers of $3$ are already forbidden, so $\exp C_T(e)\le3$.

Lemma~\ref{rankfour} applied to $E$ acting on $T$ gives $\exp T\le3$, contrary to $x\in T$ having order $9$. Hence $C_P(a_0)=1$.
\end{proof}

\begin{lemma}
\label{gapfour}
Assume~\eqref{Asetup}. If $C_D(Q)$ contains $a$ of order $3$ with $C_P(a)=1$, then every $y\in D$ of order $4$ satisfies $C_P(y)\ne1$.
\end{lemma}
\begin{proof}
Suppose that $C_P(y)=1$. Retain $C_D(Q)=Q\times E$. Since the action of $y$ on $Q$ has order at most $2$, we have $y^2\in E$, an involution. Also $a\in E$, so $E$ has even order and contains an element of order $3$.

In case (3) of Lemma~\ref{Estructure}, put $U_2=O_2(E)\ne1$. The element $a$ is fixed-point-free on $U_2$, or $6\in\om(E)$. Lemma~\ref{doublefree} gives $[U_2,P]=1$, contrary to $C_D(P)=1$.

Thus only case (4) remains: $E=U\rtimes\langle z\rangle$, where $U$ is a nonzero elementary abelian $3$-group. Here $U\trianglelefteq D$, $a\in U$, and the involution $y^2\in E$ inverts $U$.

Consider the action of $U\langle y\rangle$ on $M=P/\Phi(P)$ over a splitting field. The condition $C_P(a)=1$ gives $C_M(a)=0$. Since $5\nmid|U|$, Maschke's theorem~\cite[Theorem 1.9]{Isaacs} gives a weight decomposition. Every occurring weight $\chi$ of $U$ satisfies $\chi(a)\ne1$ and is therefore nontrivial. Since $y^2$ inverts $U$, we have $\chi^{y^2}=\chi^{-1}\ne\chi$. Thus each occurring weight has a $y$-orbit of length exactly $4$.

Lemma~\ref{norm} gives a nonzero fixed point of $y$ on $M$. Coprime lifting through $\Phi(P)$ gives $C_P(y)\ne1$, a contradiction.
\end{proof}

\begin{lemma}
\label{Aconclusion}
A group with spectrum~\eqref{spectrum} cannot satisfy the structural conditions~\eqref{Asetup}.
\end{lemma}
\begin{proof}
By Corollary~\ref{A9survive}, $B$ has no element of order $9$. Choose $x\in G$ of order $9$. Lemma~\ref{noradicalnine} gives $|xR|=3$. After conjugation, assume $a_0=xR\in D$. If $a_0$ acts nontrivially on $Q$, apply Lemma~\ref{gaplift}.

Suppose that $a_0$ centralizes $Q$. Lemma~\ref{centraldrop} gives $C_P(a_0)=1$. The group $Q\times\langle a_0\rangle$ satisfies the hypotheses of Lemma~\ref{lift21}. Lemma~\ref{orbit21} yields a normal full Sylow $5$-subgroup $\widetilde P\trianglelefteq G$ centralized by $R$.

Take $g\in G$ of order $8$. After conjugation its image in $B$ lies in $D$ and, by Lemma~\ref{noeightD}, has order $1$, $2$, or $4$. Since $R$ centralizes $\widetilde P$, this image induces the action of $g$ on $\widetilde P$. A nonidentity fixed point would yield order $40$.

If the image has order $2$ and is fixed-point-free, Lemma~\ref{inversionlemma} shows that it inverts the abelian group $P$. Faithfulness of the action of $D$ on $P$ makes this image central in $D$. It commutes with $Q\times\langle a_0\rangle$ and yields order $42$ in $D$. The trivial image is also impossible, since it fixes all of $P$.

Finally, a fixed-point-free image of order $4$ is excluded by Lemma~\ref{gapfour} with $a=a_0$. All possibilities are exhausted.
\end{proof}

\section{Completion of the proof}
\begin{proof}[Proof of the theorem]
Suppose that a solvable group $G$ satisfies~\eqref{spectrum}. By Lemma~\ref{orientations}, a Hall $\{5,7\}$-subgroup has orientation (A) or (B).

In orientation (A), Lemma~\ref{Asetupderived} gives~\eqref{Asetup}, which is excluded by Lemma~\ref{Aconclusion}.

In orientation (B), Lemma~\ref{vfour} gives $C_2^2\le G$: cyclic and generalized quaternion Sylow $2$-subgroups have been excluded by Lemmas~\ref{nocyclictwo} and~\ref{noquaterniontwo}. Lemma~\ref{sevenconclusion} excludes this orientation as well.

Both possibilities lead to a contradiction, proving the theorem.
\end{proof}

\medskip
\noindent This work was carried out with the assistance of ChatGPT (OpenAI, GPT-6).

\end{document}